\documentclass[final]{siamltex}
\usepackage{stmaryrd}
\usepackage{amssymb}
\usepackage{amsmath}
\usepackage{graphicx}
\usepackage{amsfonts}
\usepackage{color}

\allowdisplaybreaks[1]

\title{Higher order strong approximations of semilinear stochastic wave equation with additive space-time white  noise
\thanks{This work was supported by NSF of China
(No.11301550 and No.11171352) and China Postdoctoral
Science Foundation (No.2013M531798 and No.2014T70779).
}}
\author{Xiaojie Wang\thanks{School of Mathematics and Statistics, and School of Geosciences and Info-Physics,
       Central South University, Changsha 410083, Hunan, China
({\tt x.j.wang7@gmail.com, x.j.wang7@csu.edu.cn}).}
        \and Siqing Gan\thanks{School of Mathematics and
Statistics, Central South University, Changsha 410083,
Hunan, China  ({\tt sqgan@csu.edu.cn}).}
\and Jingtian Tang\thanks{School of Geosciences and Info-Physics, Central South University, Changsha 410083,
Hunan, China  ({\tt jttang@csu.edu.cn}).}
}

\begin{document}

\maketitle

\begin{abstract}
Novel fully discrete schemes are developed to numerically approximate a semilinear stochastic wave equation driven by additive space-time white noise. Spectral Galerkin method is proposed for the spatial discretization, and exponential time integrators involving linear functionals
of the noise are introduced for the temporal approximation. The resulting fully discrete schemes are very easy to implement and allow for higher strong convergence rate in time than existing time-stepping schemes such as the Crank-Nicolson-Maruyama scheme and the stochastic trigonometric method. Particularly, it is shown that the new schemes achieve in time an order of $1- \epsilon$ for arbitrarily small $\epsilon >0$, which exceeds the barrier order $\frac{1}{2}$ established by Walsh \cite{Walsh06}. Numerical results confirm higher convergence rates and computational efficiency of the new schemes.
\end{abstract}



%\pagestyle{myheadings} \thispagestyle{plain}
%\markboth{P. DUGGAN AND
%V. A. U. THORS}{SIAM MACRO EXAMPLES}

\section{Introduction}
\label{sect:intro}
Wave motions and mechanical vibrations are two common physical phenomena that are usually modeled by hyperbolic partial differential equations.
In many practical applications, random perturbation occurs and a noisy force term is hence included in
the model problems. This leads to stochastic partial differential equations (SPDEs) of hyperbolic type \cite{Cpl07,Walsh86}. One of the fundamental hyperbolic SPDEs is the stochastic wave equation, which describes a variety of physical processes, such as the motion of a vibrating string \cite{CE72} and the motion of a strand of DNA \cite{DKMNX09}.

The present work deals with the strong approximations (cf. \cite{KP92}) of the semilinear stochastic wave equation (SWE) driven by additive space-time white noise,
\begin{equation}\label{eq:SWE}
\left\{
    \begin{array}{lll}
    \frac{\partial^2 u}{\partial t^2} = \frac{\partial^2 u}{\partial x^2} + f(x, u) + \dot{W}, \quad  t \in (0, T], \:\: x \in (0,1),
    \\
     u(0, x) = u_0(x), \, \frac{\partial u}{\partial t}  (0,x) = v_0(x), \: x \in (0,1),
     \\
     u(t, 0) = u(t,1) = 0,  \: t >0,
    \end{array}\right.
\end{equation}
%where $\dot{W}$ is a noise in both time and space,
%($\mathbb{R} \times (0,1)$)
where $T \in (0, \infty)$ and $f \colon [0,1] \times \mathbb{R} \rightarrow \mathbb{R}$ is a smooth nonlinear function satisfying \begin{align}\label{f_condition1}
|f(x, z)| \leq & L(|z|+1), \\
\label{f_condition2}
\big| \tfrac{\partial f}{\partial z}(x,z)\big| \leq L,
\quad \big| \tfrac{\partial^2 f}{\partial x \partial z } (x,z)\big|&\leq L, \quad  \mbox{and} \quad  \big| \tfrac{\partial^2 f}{\partial z^2}(x,z)\big| \leq L
\end{align}
for all $x \in [0,1],\, z\in \mathbb{R}$ and some constant $L>0$. The initial data $u_0$ and $v_0$ are random variables that will be specified later.
The forcing term $\dot{W}$ is a space-time
white noise (see below), which best models the fluctuations generated
by microscopic effects in a homogeneous physical system
\cite{GC83}.
%That is, throughout this work we
%restrict ourselves to the cylindrical $I$-Wiener process.
%We refer to \cite{Cpl07,DKMNX09,DZ92,PR07,Walsh86} for
%more relevant discussion on cylindrical $I$-Wiener process.

In recent years, much progress has been made
in both strong and weak approximations of parabolic
SPDEs, see
\cite{JK09b,JK11,Kruse12,Lindgren12} and the references therein.
%for these and other developments.
In contrast to the parabolic case, there exist only a very limited number of works devoted to the numerical study of stochastic wave equations
\cite{CL07,CLS13,HE10,KLS10,KLL11,KLL12,
MPW03,QS06,SH08,Walsh06}. In \cite{QS06}, a finite difference method was considered for spatial semi-discretizations of SWE subject to multiplicative space-time white noise and a strong convergence rate of order $\tfrac{1}{3}$ was obtained. The convergence rate was improved from $\frac 1 3$ to $\frac 1 2-\epsilon$ for arbitrarily small $\epsilon>0$ in \cite{CL07} by using a spectral Galerkin method in spatial approximation of SWE with additive noise.  Using an adaptation of ``leapfrog'' discretization, Walsh \cite{Walsh06} constructed a (fully discrete) finite difference scheme,
which attains convergence order $\tfrac{1}{2}$ in both time and space. In a series of works on numerical analysis of linear stochastic evolution equations with additive noise \cite{KLS10,KLL11,KLL12}, spatial approximations by finite element methods and temporal discretizations by rational approximations to the exponential function have been investigated. For the case of space-time white noise, the strong convergence results in \cite{KLS10,KLL12} (Theorem 5.1 in \cite{KLS10} and Theorem 4.6 in \cite{KLL12}) imply convergence rates of $\tfrac{r}{r+1} \beta$ in space and $\tfrac{p}{p+1} \beta$ in time, for $\beta < \tfrac{1}{2}$ and $p,r \in \mathbb{N}^+$ being method parameters. Recently in \cite{CLS13}, a stochastic trigonometric method was introduced for the temporal approximation of linear SWEs, which strongly converges with order $\tfrac{1}{2}-\epsilon$ for arbitrary small $\epsilon >0$ in the space-time white noise case (see Theorem 4.1 in \cite{CLS13}).

%According to our earlier discussion, there exist
%no numerical schemes in the literature of strong
%convergence rates (in time and in space) exceeding
%$\tfrac{1}{2}$, which seems to be an order barrier.
To the best of our knowledge, we have not found any numerical method that strongly converges with a rate faster than $\tfrac{1}{2}$ in the literature. This seems to be an order barrier.
In fact, the limit on the convergence rate of numerical schemes for SWEs driven by space-time white noise, has been established \cite{Walsh06} in the sense that no scheme based on the basic increments of white
noise strongly converges at a rate faster than
$\tfrac{1}{2}$. An interesting question thus
arises as to whether it is possible to overcome
the order barrier.
%Unfortunately, up to now no positive result has
%been reported.
In this work, we provide a positive answer to this question by designing two fully discrete schemes for the SWE \eqref{eq:SWE}, which enjoy a strong convergence order greater than $\tfrac{1}{2}$. More precisely, we spatially discretize \eqref{eq:SWE} by a spectral Galerkin method, and then propose two exponential time integrators involving two linear functionals of the noise.
As shown in the main convergence result (Theorem \ref{thm:main.result}), under the conditions \eqref{f_condition1} and \eqref{f_condition2} the proposed fully discrete schemes strongly converge with order $\tfrac{1}{2}- \epsilon$ in space and order $1- \epsilon$ in time for arbitrarily small $\epsilon >0$. Compared with existing schemes mentioned earlier, the proposed schemes are easy to implement and produce significant improvement on the computational efficiency (see Section \ref{sect:implementation}).

Finally, we mention that the idea of
using linear functionals of the noise process in time-stepping
schemes was exploited in \cite{JK09a,JKW11}
to approximate semilinear stochastic heat equations with
additive noise.
%in contrast to the usual exponential Euler
%method involving only basic Wiener increments.
In \cite{JK09a}, a so-called accelerated exponential
Euler (AEE) method \cite{DJR10} is shown to strongly converge with order $1-\epsilon$ under seriously
restrictive commutativity assumptions
%on the derivative operator of the nonlinear operator
%$F$ defined by \eqref{eq:Nemytskij} and the Laplacian
(see Assumption 2.4 in \cite{JK09a} and discussion in the
introduction of \cite{JKW11}). In the present work, such AEE method is successfully adapted to solve the SWE \eqref{eq:SWE} and under standard assumptions the strong
convergence rate of order $1-\epsilon$ in time is proved
in the space-time white noise case.
%with the restrictive assumptions removed.
%Although the
%commutativity conditions are fulfilled for the linear
%function $f(u)=cu$, they are quite restrictive and
%exclude most nonlinear functions such as
%$f(u) = \tfrac{1+u}{1+ u^2}$.
%In this article we prove the rate $1-\epsilon$ in time of
%the schemes on conditions \eqref{f_condition1}
%and \eqref{f_condition2}, without imposing the
%seriously restrictive
%commutativity assumptions as required in \cite{JK09a}.

%When the driven noise is a smoother $Q$-Wiener process
%with $Q$ being trace class (i.e., Tr$(Q) < \infty$),
%Jentzen et al. \cite{JKW11} get rid of such commutativity
%conditions for the accelerated exponential Euler scheme.
%For the most interesting but tough case when the driven
%noise is white in both time and space, however, it seems
%extremely difficult to recover the strong convergence rate
%$1-\epsilon$ in time under standard assumptions.
%In this article we shall make attempts in this direction
%for stochastic wave equation.  We first devise two
%accelerated exponential Euler methods for SWE \eqref{eq:SWE}
%with noise being white in time as well as space ($Q=I$).
%Last, but most importantly, we prove the rate $1-\epsilon$
%in time of the schemes on conditions \eqref{f_condition1}
%and \eqref{f_condition2}, without imposing the restrictive
%commutativity assumptions as required in \cite{JK09a}.

The rest of this paper is organized as follows. In the next section, some preliminaries are collected and an abstract formulation of \eqref{eq:SWE} is set forth. In Section \ref{sect:Galerkin}, we analyze the strong approximation error arising from the spatial discretization by a spectral Galerkin method. Then two exponential time integrators are introduced and strong convergence of the fully discrete approximations are studied in Section \ref{sect:superconvergence}. Numerical experiments confirming our theoretical results are presented in Section \ref{sect:implementation}. The paper is concluded with some brief remarks in Section 6.

\section{Preliminaries and framework}
\label{sect:framework}

Let $(U,\, \langle \cdot, \cdot \rangle, \,\|\cdot \|)$
and $(H,\, ( \cdot, \cdot), \, \||\cdot \||)$
be two separable Hilbert spaces.
By $\mathcal{L}(U,H)$ we denote the Banach space of bounded linear operators from $U$ to $H$ and for short we write $\mathcal{L}(U) :=\mathcal{L}(U,U)$.
Additionally, we need the Banach space of
Hilbert-Schmidt operators, denoted by $\mathcal{L}_2(U,H)$, equipped with the norm
\begin{equation}\label{eq:normdef} \|\Gamma\|_{\mathcal{L}_2(U,H)} = \Big(\sum_{i=1}^{\infty}\||\Gamma \eta_i\||^2\Big)^{1/2},
\end{equation}
where $\{\eta_i\}_{i \in \mathbb{N}}$ is an orthonormal basis of $U$ and the norm does not depend on the particular choice of the basis \cite{DZ92,PR07}. For brevity, we write $\mathcal{L}_2(U) := \mathcal{L}_2(U,U)$. If $\Gamma_1\in \mathcal{L}(U,H)$ and $\Gamma_2\in \mathcal{L}_2(U)$, then $\Gamma_1 \Gamma_2  \in \mathcal{L}_2(U,H)$, and
\begin{equation}\label{eq:LHS}
\| \Gamma_1 \Gamma_2\|_{\mathcal{L}_2(U,H)} \leq \|\Gamma_1\|_{\mathcal{L}(U,H)} \cdot \|\Gamma_2\|_{\mathcal{L}_2(U)}.
\end{equation}
Let $\left(\Omega,\mathcal {F},\mathbb{P}\right)$ be a probability space with a normal filtration $\{\mathcal{F}_t\}_{0\leq t\leq T}$ and by $L^p(\Omega, U)$ we denote the space of $U$-valued integrable random variables with the norm defined by
$
\|\varphi\|_{L^p(\Omega, U)} = \big(\mathbb{E} \big[ \| \varphi \|^p \big] \big)^{\frac{1}{p}} < \infty
$
for any $p\geq 2$.

Next, we take $U := L^2\big((0,1), \mathbb{R}\big)$
to denote the space of real-valued square integrable
functions, equipped with the usual norm  $\| \cdot \|$
and inner product $\langle \cdot, \cdot \rangle$.
Let $- \Lambda = \Delta : \mathcal{D}(\Lambda)\subset U
\rightarrow U $ be the Laplace operator with
$\mathcal{D}(\Lambda) = H^2(0, 1) \cap H_0^1(0, 1)$, where $H^m(0,1)$ denote the standard Sobolev spaces of integer order $m \geq 1$ and $H_0^1(0, 1) := \{ \varphi \in H^1(0,1) \colon \varphi(0) = \varphi(1) =0 \}$.
Then $\Lambda$ is a densely defined, self-adjoint, positive operator with compact inverse. Moreover, the eigenvalue problem
\begin{equation}\label{eq:lambda.eigen.prob}
\Lambda e_i = \lambda_i e_i, \, i \in \mathbb{N}
\end{equation}
provides an orthonormal
basis $\{e_i=\sqrt{2}\sin(i \pi x), \,x
\in(0,1)\}_{i \in \mathbb{N}}$ for $U$ and an increasing sequence of eigenvalues $\lambda_i=\pi^2 i^2, i \in \mathbb{N}$.
Additionally, let $F: U \rightarrow U$ be a Nemytskij operator associated to $f$ as in \eqref{eq:SWE}, defined by
\begin{equation}\label{eq:Nemytskij}
F (\varphi)(x) = f(x, \varphi(x)), \quad x \in (0,1), \: \varphi \in U.
\end{equation}
Now one can rewrite \eqref{eq:SWE} as an
abstract form in It\^{o}'s sense
%\textcolor{red}{@@ u(t) or u(x)? @@}
\begin{align}\label{eq:Abstr.SWE}
\left\{\!
    \begin{array}{ll}
    \mbox{d} \dot{u} = -\Lambda u \mbox{d}t + F(u) \mbox{d}t + \mbox{d}W(t), \quad t \in (0, T],
    \\
     u(0) = u_0, \, \dot{u}(0) = v_0,
    \end{array}\right.
\end{align}
where $u$ is regarded as a $U$-valued stochastic process and $\dot{u}$ stands for the time derivative of $u$. The driven stochastic process $W(t)$ is a cylindrical $I$-Wiener process with respect to $\{\mathcal
{F}_t\}_{0\leq t\leq T}$,
which can be represented as follows \cite{DZ92,PR07}:
\begin{align}\label{W.representation}
W(t) = \sum_{i=1}^{\infty} \beta_i(t) e_i, \quad t \in [0, T],
\end{align}
where $\{\beta_i(t) \}_{i\in \mathbb{N}}$ are independent real-valued Brownian motions and $\{ e_i \}_{i \in \mathbb{N}}$ are the eigenvectors of $\Lambda$ defined by \eqref{eq:lambda.eigen.prob}.
Moreover, note that the derivative operators of $F$ are given by
\begin{align} \label{eq:F.Deriv}
F'(\varphi)(\psi) \,(x)= &  \tfrac{\partial f}{\partial z}(x,\varphi(x)) \cdot \psi(x), \quad x \in (0,1),
\\
F''(\varphi)(\psi_1, \psi_2) \,(x)= &  \tfrac{\partial^2 f}{\partial z^2}(x,\varphi(x)) \cdot \psi_1(x) \cdot \psi_2(x), \quad x \in (0,1),
\label{eq:F.Deriv2}
\end{align}
for all $\varphi, \psi, \psi_1, \psi_2 \in U$. It is worthwhile to keep in mind that the derivative operators $F'(\varphi), F''(\varphi), \varphi \in U$ defined in the above way are self-adjoint. Thanks to \eqref{f_condition1} and \eqref{f_condition2}, the Nemytskij operator $F$ satisfies
\begin{align}\label{eq:Nf.condition1}
  \|F(\varphi) \| \leq& \sqrt{2} L(\|\varphi\|+1),
  \\
  \label{eq:Nf.condition2}
  \|F(\varphi_1)-F(\varphi_2) \| \leq& L \|\varphi_1-\varphi_2\|
\end{align}
for all $\varphi,\varphi_1, \varphi_2 \in U$. Also, it is obvious that
\begin{align}\label{eq:Lambda.Q}
\big\|\Lambda^{\frac{\beta-1}{2}} \big\|^2_{\mathcal{L}_2(U)} = \sum_{i=1}^{\infty} \| \Lambda^{\frac{\beta-1}{2}}e_i \|^2 =  \sum_{i=1}^{\infty} \frac{\pi^{2(\beta-1)}} {i^{2(1-\beta)}}  \leq \bar{c}_{\beta} < \infty, \quad \mbox{for any } \beta <\tfrac{1}{2}.
\end{align}

%Here and throughout this work, $C$ appearing in the
%following estimates is a generic constant that may
%vary from one place to another and depends only on
%$T, \epsilon, L$ and initial data, but is independent
%of \textcolor{red}{$\tau, N$}.

%the existence and uniqueness of the mild solution of
%the stochastic wave equation \eqref{eq:Abstr.SWE}

In order to define the mild solution of \eqref{eq:Abstr.SWE} appropriately, we shall reformulate  \eqref{eq:Abstr.SWE} as a stochastic evolution equation in a new Hilbert space $H$ to fall into the semigroup framework in \cite{DZ92}. To this end, we need  additional spaces and notations. The above setting enables us to define fractional powers of $\Lambda$ in a simple way (see, e.g., \cite[Appendix B.2]{Kruse12}). Accordingly, we introduce the separable Hilbert space $\dot{H}^{\alpha} := \mathcal{D}(\Lambda^{\frac{\alpha}{2}})$, $\alpha \in \mathbb{R}$, equipped with the inner product
\begin{equation}\label{eq:Lambda.norm}
\langle \varphi, \psi \rangle_{\alpha} := \langle \Lambda^{\frac{\alpha}{2}} \varphi, \Lambda^{\frac{\alpha}{2}} \psi\rangle =  \sum_{i=1}^{\infty} \lambda_i^{\alpha} \langle \varphi, e_i\rangle \langle \psi, e_i\rangle,
\quad
\varphi, \psi \in \dot{H}^{\alpha},
\end{equation}
where $\{ (\lambda_i, e_i) \}_{i=1}^{\infty}$ are the
eigenpairs of $\Lambda$. The corresponding norm is defined by $\|\varphi\|_{\alpha} = \sqrt{\langle \varphi, \varphi \rangle _{\alpha}}$ for $\varphi \in \dot{H}^\alpha$.
Then $\dot{H}^{0} = U$ and $\dot{H}^{\alpha} \subset \dot{H}^{\beta}$ if $\alpha \geq \beta$. Moreover, $\dot{H}^{-\gamma}$ can be identified with the dual space $\big(\dot{H}^{\gamma}\big)^*$ for $\gamma >0$ \cite{TV97}.  Further, we introduce the product space
$H^{\alpha} := \dot{H}^{\alpha} \times \dot{H}^{\alpha-1}$, $\alpha \in \mathbb{R}$,
endowed with the inner product
\begin{equation}\label{eq:prod.inner.prod}
( Y, \hat{Y} )_{\alpha} := \langle \varphi,   \hat{\varphi} \rangle_{\alpha} + \langle \psi, \hat{\psi} \rangle_{\alpha-1}, \quad Y = (\varphi, \psi)^T, \: \hat{Y} = (\hat{\varphi}, \hat{\psi})^T,
\end{equation}
and the usual norm
\begin{align}\label{eq:product.norm}
\|| Y \||_{\alpha}^2 := \|\varphi\|_{\alpha}^2 + \| \psi\|^2_{\alpha-1}, \quad Y=(\varphi,\psi)^T.
\end{align}
It is easy to check that $(H^{\alpha},\, ( \cdot, \cdot)_{\alpha} )$, $\alpha \in \mathbb{R}$,
is a separable Hilbert space.
For the special case $\alpha =0$, we denote $H := H^0 = \dot{H}^0 \times \dot{H}^{-1}$,
$(\cdot, \cdot) := (\cdot, \cdot)_{0}$, and
$\|| \cdot \|| := \|| \cdot \||_0$.
%and define the inner product in $H$ by $(Y, \hat{Y}) =
%\langle \varphi, \hat{\varphi} \rangle + \langle
%\Lambda^{-\frac{1}{2}} \psi, \Lambda^{-\frac{1}{2}}
%\hat{\psi} \rangle $ for $Y = (\varphi, \psi)^T$,
%$\hat{Y} = (\hat{\varphi}, \hat{\psi})^T$.

At this point, we introduce the velocity of the solution $u$, denoted by $v = \dot{u}$, and formally transform \eqref{eq:Abstr.SWE} into the following Cauchy problem
\begin{equation}\label{eq:SEE}
\left\{
    \begin{array}{ll}
    \mbox{d} X(t) = A X(t) \mbox{d}t + \mathbf{F}(X) \mbox{d}t + B\mbox{d}W(t), \quad  t \in (0, T],
    \\
     X(0)= X_0,
    \end{array}\right.
\end{equation}
where $X_0 = (u_0,v_0)^T$ and
\begin{equation}
X = \bigg[ \!\begin{array}{c}
    u
    \\
    v
    \end{array}\!\bigg], \:
    A = \bigg[\! \begin{array}{cc}
    0 & I
    \\
    -\Lambda & 0
    \end{array}\!\bigg],\:
    \mathbf{F}(X) = \bigg[\! \begin{array}{c}
    0
    \\
    F(u)
    \end{array}\!\bigg],\:
    B = \bigg[\! \begin{array}{c}
    0
    \\
    I
    \end{array}\!\bigg].
\end{equation}
Here $X_0$ is assumed to be an $\mathcal {F}_0$-measurable $H$-valued random variable and $B$ is considered as an operator from $\dot{H}^{-1}$ to $H$. From now on, we regard $\Lambda$ as an operator from $\dot{H}^1$ to $\dot{H}^{-1}$, defined by $( \Lambda \varphi )(\psi) = \langle \nabla \varphi, \nabla \psi \rangle$ for $\varphi, \psi \in \dot{H}^1$ and define the domain of $A$ by
\begin{equation*}\label{eq:Domain.A}
  \mathcal{D} (A) = \bigg\{ Y = (\varphi, \psi)^T
   \in H: A Y =  \bigg[\! \begin{array}{c}
    \psi
    \\
    -\Lambda \varphi
    \end{array}\!\bigg] \in H
    = \dot{H}^0 \times \dot{H}^{-1}
    \bigg\}
    = H^1 = \dot{H}^1 \times \dot{H}^{0}.
\end{equation*}
Then the operator $A$ is the generator of a strongly continuous semigroup $E(t), t \geq 0$ on $H$ \cite[Section 5.3]{Lindgren12}, that can be written as
\begin{align}\label{eq:Et}
E(t) = e^{tA} = \left[\! \begin{array}{cc}
    C(t) & \Lambda^{-\frac{1}{2}} S(t)
    \\
    -\Lambda^{\frac{1}{2}}S(t) & C(t)
    \end{array} \!\right].
\end{align}
Here $C(t)\!=\! \cos(t\Lambda^{\frac{1}{2}})$ and
$S(t)\!=\! \sin(t\Lambda^{\frac{1}{2}})$ are the so-called
cosine and sine operators, which can be expressed in
terms of the eigenpairs $\{ \lambda_i, e_i \}_{i \in
\mathbb{N} }$:
%For example, we have
\begin{align*}
C(t) \varphi =\sum_{i=1}^{\infty} \cos(t \lambda_i^{1/2}) \langle \varphi, e_i \rangle e_i, \quad  S(t) \varphi =\sum_{i=1}^{\infty}  \sin(t \lambda_i^{1/2}) \langle \varphi, e_i \rangle e_i
\end{align*}
for $t \geq0, \varphi \in \dot{H}^{-1}$.
%Here $( \cdot, \cdot)$ denotes the pairing between
%$\dot{H}^{-1}$ and $\dot{H}^{1}$.
Before proceeding further, we briefly state some properties of $C(t)$ and $S(t)$, which will be used frequently later. As defined above, the cosine and sine operators are bounded in the sense that $\| C(t) \varphi \| \leq \| \varphi \|$ and $\| S(t) \varphi \| \leq \| \varphi\|$ hold for all $\varphi \in U$. In addition, these two operators satisfy the trigonometric identity $\|S(t)\varphi\|^2 + \|C(t)\varphi\|^2 = \|\varphi\|^2$ for $\varphi \in U$. Moreover, $\Lambda^{\gamma}$, $\gamma \in \mathbb{R}$ commutes with $C(t), S(t)$ since they are all defined in terms of the eigenbasis $\{e_i\}_{i \in \mathbb{N}}$. With the aid of these properties together, one can show that $\|| E(t) Y \||_{\gamma} \leq \|| Y \||_{\gamma}$  for $t \in [0, \infty), \gamma \in \mathbb{R},  Y \in H^{\gamma}$. Now, we look at the existence and uniqueness of the mild solution of \eqref{eq:SEE}, which has been discussed in \cite{CN88,QS06} using different frameworks.

\begin{theorem}\label{thm:unique.mild}
Suppose conditions \eqref{f_condition1} and \eqref{f_condition2} are fulfilled, let $W(t), t\in [0, T]$ be the cylindrical $I$-Wiener process represented by \eqref{W.representation}, and let $X_0$ be an $\mathcal {F}_0$-measurable $H$-valued random variable satisfying
$
\|X_0\|_{L^p(\Omega, H)} < \infty
$
%\begin{equation}
%\begin{array}{l}
%\|X_0\|_{L^p(\Omega, H^{\frac{1}{2}})} =
%\Big( \mathbb{E} \|| X_0\||_{\frac{1}{2}}^p
%\Big)^{\frac{1}{p}} = \left( \mathbb{E} \Big[\!
%\left( \| u_0\|_{\frac{1}{2}}^2 + \| v_0\|_{-\frac{1}{2}}^2
%\right)^{p/2}\Big]\right)^{\frac{1}{p}} < \infty
%\end{array}
%\end{equation}
for some $p \geq2$. Then SWE \eqref{eq:SEE} has a unique mild solution given by \begin{align}\label{eq:mild.solution}
\small
  X(t) = E(t)X_0 +  \int_0^t E(t-s) \mathbf{F}(X(s)) \, \mbox{d}s + \int_0^t E(t-s)B \, \mbox{d} W(s) \quad a.s.
\end{align}
for each $t \in [0, T]$.
Additionally, if $
\|X_0\|_{L^p(\Omega, H^{1/2})} < \infty
$ then there exists a constant $K_1(p,\beta, T) \in [0, \infty)$ depending on $p, \beta, T$ such that for any $0\leq \beta < \frac{1}{2}$
\begin{align}\label{eq:X.MB}
  \sup_{t \in [0, T]} \|X(t)\|_{L^p(\Omega, H^{\beta})} \leq K_1 \big( \|X_0\|_{L^p(\Omega, H^{\beta})} + 1 \big).
\end{align}
\end{theorem}
{\it Proof.} We first claim that the nonlinear operator $\mathbf{F}\colon H \rightarrow H$ satisfies the globally Lipschitz condition and linear growth condition.
In fact, for any $Y=(\varphi, \psi)^T \in H$, $Y_1=(\varphi_1, \psi_1)^T \in H$,
$Y_2=(\varphi_2, \psi_2)^T \in H$, we infer
\begin{align*}%\label{}
\|| \mathbf{F}(Y_1)-\mathbf{F}(Y_2) \||
= & \|\Lambda^{-\frac{1}{2}} \big(F(\varphi_1)- F(\varphi_2)\big)\|\leq L \|\varphi_1-\varphi_2\| \leq L\||Y_1-Y_2\||,
\\
\|| \mathbf{F}(Y)\||
= &
\|\Lambda^{-\frac{1}{2}} F(\varphi)\|\leq \sqrt{2}L (\|\varphi\|+1) \leq \sqrt{2}L(\||Y\||+1)
\end{align*}
by \eqref{eq:Nf.condition1}, \eqref{eq:Nf.condition2}, the definitions of $\mathbf{F}$ and $\||\cdot \||$, and stability of $\Lambda^{-\frac{1}{2}}$. Here and below by ``stability'' of a linear operator $\Gamma \colon U \rightarrow U$, we mean $\Gamma$ is bounded and satisfies $\| \Gamma \varphi \| \leq \|\varphi\|$ for all $\varphi \in U$.  Furthermore, we have
\begin{equation}
\|B\|_{\mathcal{L}_2(U,H)} = \Big( \sum_{i=1}^{\infty} \||B e_i\||^2 \Big)^{\frac{1}{2}}
=
\Big( \sum_{i=1}^{\infty} \|\Lambda^{-\frac{1}{2}} e_i\|^2 \Big)^{\frac{1}{2}}
=
\Big( \sum_{i=1}^{\infty} \lambda_i^{-1} \Big)^{\frac{1}{2}}
<
\infty.
\end{equation}
In view of Theorem 7.4 from \cite{DZ92}, one can derive the existence and uniqueness of the mild solution \eqref{eq:mild.solution}, which satisfies
\begin{align}\label{eq:X.MB2}
\sup_{t \in [0, T]} \| X(t) \|_{L^p(\Omega, H)} \leq C_{p,T} \left( \| X_0 \|_{L^p(\Omega, H)} + 1 \right).
\end{align}
By using the trigonometric identity and \eqref{eq:Lambda.Q} one can easily show that, for $ t \in [0, T]$
\begin{equation}\label{eq:stoch.convol.MB}
\int_0^t \|E(t-s) B\|^2_{\mathcal{L}_2(U,H^{\beta})} \mbox{d}s \leq t \|\Lambda^{\frac{\beta-1}{2}}\|^2_{\mathcal{L}_2(U)} \leq t \bar{c}_{\beta}  < \infty, \quad \beta \in [0, \tfrac{1}{2}).
\end{equation}
Using  this together with \eqref{eq:Nf.condition1}, \eqref{eq:product.norm}, the Burkholder-Davis-Gundy type inequality (\cite[Lemma 7.2]{DZ92}) and some properties of the operators $\Lambda^{\gamma}$, $\gamma \in \mathbb{R}$, $C(t)$ and $S(t)$ mentioned earlier, we deduce from \eqref{eq:mild.solution} that
\begin{align}
&\|X(t)\|_{L^p(\Omega, H^{\beta})}
\nonumber \\
\leq &
\|E(t)X_0\|_{L^p(\Omega, H^{\beta})} +  \int_0^t \left\|E(t-s) \mathbf{F}(X(s))\right\|_{L^p(\Omega, H^{\beta})} \mbox{d}s
\nonumber \\
& +
\left\|\int_0^t E(t-s)B \, \mbox{d} W(s)\right\|_{L^p(\Omega, H^{\beta})}
\nonumber \\
\leq &
\big\| \||E(t)X_0\||_{\beta} \big\|_{L^p(\Omega, \mathbb{R})} \!+
c_p \bigg(\! \int_0^t \| E(t-s) B \|^2_{\mathcal{L}_2(U, H^{\beta})} \mbox{d}s \! \bigg)^{\frac{1}{2}}
\nonumber \\
& +
\int_0^t \Big\| \Big( \| \Lambda^{-\frac{1}{2}} S(t-s)F(u(s))\|^2_{\beta} + \| C(t-s)F(u(s))\|^2_{\beta-1}\Big)^{\frac{1}{2} }  \Big\|_{L^p(\Omega, \mathbb{R})} \mbox{d}s
\nonumber \\ \leq&
\big\| \||X_0\||_{\beta} \big\|_{L^p(\Omega, \mathbb{R})} + \sqrt{t \bar{c}_{\beta} }\, c_p
+ \int_0^t \big\| \, \| F(u(s)) \|_{\beta-1} \big\|_{L^p(\Omega, \mathbb{R})} \mbox{d}s
\nonumber \\ \leq&
\|X_0\|_{L^p(\Omega, H^{\beta})} + \sqrt{t \bar{c}_{\beta} }\, c_p + \sqrt{2} L \int_0^t \big\| \, \| u(s)\| +1 \big\|_{L^p(\Omega, \mathbb{R})} \mbox{d}s
\nonumber \\
\leq&
\|X_0\|_{L^p(\Omega, H^{\beta})} + \sqrt{t \bar{c}_{\beta} }\, c_p  + \sqrt{2} L t + \sqrt{2} L \int_0^t \|u(s)\|_{L^p(\Omega, U)} \mbox{d}s
\label{eq:X.MB4}
\end{align}
for $ t \in [0, T]$ and $0\leq \beta <\frac{1}{2}$.
%and where the stability properties of the operators
%$S(t), C(t), \Lambda^{-\gamma}$, $\gamma \in [0, \infty)$
%in $U$ and the fact that $\Lambda^{-\gamma}$, $\gamma
%\geq 0$ commutes with $C(t), S(t)$ were also used.
The definition of the norm $\|| \cdot \||$ and \eqref{eq:X.MB2} guarantee
\begin{equation}
\|u(t)\|_{L^p(\Omega, U)} \leq  \| X(t) \|_{L^p(\Omega, H)} \leq C_{p,T} \left( \| X_0 \|_{L^p(\Omega, H)} + 1 \right)
\end{equation}
for all $t \in [0, T]$. This and \eqref{eq:X.MB4} together thus yield the desired estimate \eqref{eq:X.MB}.
$\square$

Replacing $E(t)$ by \eqref{eq:Et} one can write \eqref{eq:mild.solution} as
\begin{align}\label{eq:mild.concrete}
\left\{\!
    \begin{array}{l}
    u(t)= C(t) u_0 + \Lambda^{-\frac{1}{2}} S(t)v_0 + \int_0^t \Lambda^{-\frac{1}{2}} S(t-s)F(u(s)) \mbox{d}s + \mathcal{O}_t,
    \\
    v(t)= -\Lambda^{\frac{1}{2}} S(t) u_0 + C(t)v_0 + \int_0^t C(t-s)F(u(s)) \mbox{d}s +  \widehat{\mathcal{O}}_t,
    \end{array}\right.
\end{align}
where $t \in [0, T]$ and we used the notations
\begin{equation}\label{eq:stoch.convol.notat}
\begin{array}{l}
\mathcal{O}_t = \int_0^t \Lambda^{-\frac{1}{2}} S(t-s)\, \mbox{d} W(s), \quad \widehat{\mathcal{O}}_t = \int_0^t C(t-s) \,\mbox{d} W(s).
\end{array}
\end{equation}

% Next we are to design approximations of the mild
% solution $u(t)$ and measure the discrepancy between
% the approximations and $u(t)$. To this end, one has to
% discretize both the time interval $[0,T]$ and the
% infinite dimensional space $U$.
% For temporal discretizations of SWE, the finite
% difference method is a common choice
% \cite{HE10,KLL12,MPW03,Walsh06}, while
% spatial discretizations can be achieved with finite
% difference\cite{HE10,QS06,Walsh06}, finite element
% \cite{KLS10,KLL11,KLL12} and spectral Galerkin
% \cite{CL07,SH08} methods. In this work we consider
% finite difference time discretization and spectral
% Galerkin discretization in space.

%In terms of ...

\section{The spectral Galerkin approximation of SWE}
\label{sect:Galerkin}

In this section we consider the spatial discretizations of \eqref{eq:SEE}  by a spectral Galerkin method.
To this end, for $N\in \mathbb{N}$ we define a finite dimensional subspace of $U$ by
$
%\begin{equation}\label{eq:space.HN}
 U_N := \mbox{span} \{e_1, e_2, \cdots, e_N \},
%\end{equation}
$
and the projection operator $P_N \colon \dot{H}^{\alpha}\rightarrow U_N$ by
\begin{align}
  P_N \xi = \sum_{i=1}^N \langle \xi, e_i \rangle e_i, \quad \forall\, \xi \in \dot{H}^{\alpha}, \, \alpha \geq -1.
\end{align}
The definition of $P_N$ immediately implies
\begin{equation}
\| P_N \varphi \|^2 = \Big\| \sum_{i = 1}^N \langle \varphi, e_i\rangle e_i \Big\|^2 = \sum_{i = 1}^N | \langle \varphi, e_i\rangle |^2 \leq \sum_{i = 1}^{\infty} | \langle \varphi, e_i\rangle |^2 = \|\varphi\|^2, \quad \forall \: \varphi \in U.
\end{equation}
We emphasize that $U_N$ here is chosen as the linear space spanned by the $N$ first eigenvectors of $\Lambda$. This ensures easy simulations of the stochastic convolutions in the proposed numerical schemes (see below). Now we define $\Lambda_N \colon U_N\rightarrow U_N$ by
\begin{align}\label{eq:Lambda.Galerkin}
\Lambda_N \xi = \Lambda P_N \, \xi =P_N \Lambda  \, \xi = \sum_{i=1}^N \lambda_i \langle \xi, e_i \rangle e_i, \quad \forall \: \xi \in U_N.
\end{align}
Similarly, one can define $\Lambda_N^{\gamma} \colon U_N \rightarrow U_N $,
$\gamma \in \mathbb{R}$ in $U_N$ as
$
\Lambda_N^{\gamma} \xi : = \sum_{i=1}^N \lambda_i^{\gamma} \langle \xi, e_i \rangle e_i,  \: \xi \in U_N.
$
We next apply the spectral Galerkin method to \eqref{eq:SEE}. This gives finite dimensional stochastic differential equations (SDEs) in $H_N:=U_N\times U_N$
\begin{align}\label{eq:Galerkin.SEE}
  \left\{
    \begin{array}{ll}
    \mbox{d} X^N(t) = A_N X^N(t) \mbox{d}t + \mathbf{F}_N(X^N) \mbox{d}t + B_N\mbox{d}W(t), \quad t \in (0, T],
    \\
     X^N(0)= X^N_0,
    \end{array}\right.
\end{align}
where $X^N_0 = (P_N u_0, P_N v_0)^T$ and
\begin{align*}
X^N = \bigg[\! \begin{array}{c}
    u^N
    \\
    v^N
    \end{array} \!\bigg], \:
    A_N = \bigg[\! \begin{array}{cc}
    0 & I
    \\
    -\Lambda_N & 0
    \end{array} \!\bigg], \:
    \mathbf{F}_N(X^N)=\bigg[\! \begin{array}{c}
    0
    \\
    P_N\,F(u^N)
    \end{array} \!\bigg],\:
    B_N=\bigg[\! \begin{array}{c}
    0
    \\
    P_N
    \end{array}\!\bigg].
\end{align*}
Analogously, the operator $A_N$ is the generator of a strongly continuous semigroup $E_N(t), t \geq 0$ on $U_N \times U_N $ and
\begin{align}\label{eq:ENt}
E_N(t) = e^{tA_N} \!=\! \bigg[\! \begin{array}{cc}
    C_N(t) & \Lambda_N^{-\frac{1}{2}} S_N(t)
    \\
    -\Lambda_N^{\frac{1}{2}}S_N(t) & C_N(t)
    \end{array}\!\bigg],
\end{align}
where $C_N(t) \!=\! \cos(t\Lambda_N^{\frac{1}{2}} )$ and $S_N(t) \!=\! \sin(t\Lambda_N^{\frac{1}{2}} )$ for $t \geq 0$ are the cosine and sine operators defined in $U_N$. It can be verified straightforwardly that
\begin{align}\label{eq:AS.Comm}
C_N(t) P_N \varphi = C(t) P_N \varphi = P_N C(t)  \varphi,
\quad
S_N(t) P_N \varphi = S(t) P_N \varphi = P_N S(t)  \varphi
\end{align}
for $\varphi \in \dot{H}^{\alpha}, \, \alpha \geq -1$.
The following result ensures a unique global solution  of \eqref{eq:Galerkin.SEE}.

\begin{theorem}\label{thm:Galerkin.unique.mild}
Assume that all conditions in Theorem \ref{thm:unique.mild} are fulfilled. Then \eqref{eq:Galerkin.SEE} has a unique solution given by
\begin{equation}\label{eq:Galerk.mild.solut}
\small
  X^N(t) = E_N(t)X^N_0 +  \int_0^t E_N(t-s) \mathbf{F}_N(X^N(s)) \mbox{d}s + \int_0^t E_N(t-s)B_N \mbox{d} W(s) \quad a.s.
\end{equation}
for any $t \in [0, T]$. Additionally, there exists a constant $K_2(\beta, p, T)$ depending on $ \beta, p, T$ such that for any $0\leq \beta < \frac{1}{2}$ and $t \in [0, T]$,
\begin{align}\label{eq:XN.MB}
  \|X^N(t)\|_{L^p(\Omega, H^{\beta})} \leq K_2 \big( \|X_0\|_{L^p(\Omega, H^{\beta})} + 1 \big).
\end{align}
\end{theorem}
The proof of Theorem \ref{thm:Galerkin.unique.mild} goes along the same lines as that of Theorem \ref{thm:unique.mild} and is thus omitted here.

Similarly to \eqref{eq:mild.concrete}, \eqref{eq:Galerk.mild.solut} can be rewritten as
\begin{align}\label{eq:Galerkin.mild.concrete}
\!\! \left\{ \!\!
    \begin{array}{l}
    u^N(t)= C_N(t) u_0^N + \Lambda_N^{-\frac{1}{2}} S_N(t)v_0^N \!+ \!\int_0^t \Lambda_N^{-\frac{1}{2}} S_N(t-s)P_N F(u^N(s)) \mbox{d}s \! + \!\mathcal{O}^N_t\!,
    \\
    v^N(t)= -\Lambda_N^{\frac{1}{2}} S_N(t) u_0^N + C_N(t)v_0^N \!+ \!\int_0^t C_N(t-s)P_N F(u^N(s)) \mbox{d}s \! + \!\widehat{\mathcal{O}}^N_t\!,
    \end{array}
    \right.
\end{align}
where for simplicity of presentation we denote $u_0^N = P_N u_0$, $v_0^N = P_N v_0$ and
\begin{equation}\label{eq:stoch.convol.Galerkin}
\begin{array}{l}
\mathcal{O}^N_t = \int_0^t \Lambda_N^{-\frac{1}{2}} S_N(t-s)P_N \, \mbox{d} W(s), \quad  \quad \widehat{\mathcal{O}}^N_t = \int_0^t C_N(t-s)P_N \, \mbox{d} W(s).
\end{array}
\end{equation}
%Before going to the spatial error analysis, we present
The following lemma is an immediate consequence of \eqref{eq:Nf.condition1}, \eqref{eq:X.MB} and \eqref{eq:XN.MB}.
\begin{lemma} \label{lemma:MB.FXFXN}
Assume that all conditions in Theorem \ref{thm:unique.mild} are fulfilled, and let $u(t)$ and $u^N(t)$ be given by \eqref{eq:mild.concrete} and \eqref{eq:Galerkin.mild.concrete}, respectively. Then there exists a constant $K_3(p, T, L, \beta)$ depending on $p, T, L, \beta$ such that
\begin{align}\label{eq:MB.FXFXN}
\big \|F (u(t))\big\|_{L^p(\Omega, U)} +  \big \|F (u^N(t))\big\|_{L^p(\Omega, U)}  \leq K_3 \big( \|X_0\|_{L^p(\Omega, H)} + 1 \big), \quad t \in [0, T].
\end{align}
\end{lemma}
Armed with the above preparations, we are now able to analyze the spatial discretization error.

\begin{theorem}[Spatial discretization error] \label{thm:spatial.error}
Suppose that all conditions in Theorem \ref{thm:unique.mild} are satisfied. Then it holds for all $t \in [0, T]$ that
\begin{align}\label{eq:u.spatial.error}
\|u^N(t)- u(t)\|_{L^2(\Omega,U)} \leq& K_4 \big(\|X_0\|_{L^2(\Omega, H^{1/2-\epsilon})} +1 \big)N^{-\frac{1}{2}+\epsilon}
\end{align}
for arbitrarily small $\epsilon >0$,
where $u(t)$ and $u^N(t)$ are given by \eqref{eq:mild.concrete} and \eqref{eq:Galerkin.mild.concrete}, respectively, and where $K_4(\epsilon, T)  \in [0, \infty)$ is a constant depending on $\epsilon, T$.
\end{theorem}

{\it Proof.}
The definitions of $u(t)$ in \eqref{eq:mild.concrete}
and $u^N(t)$ in \eqref{eq:Galerkin.mild.concrete} yield,
%\begin{align}\label{eq:Diff.Space1}
%u^N(t) - u(t) =& \big(C_N(t) P_N - C(t)\big) u_0
%+ \big( \Lambda_N^{-\frac{1}{2}}S_N(t)P_N -
%\Lambda^{-\frac{1}{2}}S(t) \big) v_0
%\nonumber \\
%& + \smallint_0^t \big(\Lambda_N^{-\frac{1}{2}}
%S_N(t-s)P_N F(u^N(s))-\Lambda^{-\frac{1}{2}}S(t-s)
%F(u(s))  \big) \mbox{d}s
%\nonumber \\
%& +  \smallint_0^t \big(\Lambda_N^{-\frac{1}{2}}
%S_N(t-s)P_N - \Lambda^{-\frac{1}{2}} S(t-s)\big)
%\mbox{d} W(s)
%\end{align}
%and
%\begin{align}\label{eq:Diff.Space2}
%v^N(t) - v(t) =& \left(-\Lambda_N^{\frac{1}{2}}
%S_N(t) P_N + \Lambda^{\frac{1}{2}}S(t)\right) u_0 +
%\left( C_N(t)P_N-C(t) \right) v_0
%\nonumber \\
%& + \int_0^t \left(C_N(t-s)P_N F(u^N(s)) -
%C(t-s) F(u(s)) \right) \mbox{d}s
%\nonumber \\
%& +  \int_0^t \left(C_N(t-s)P_N  - C(t-s)\right)
%\mbox{d} W(s).
%\end{align}
for all $t \in [0, T]$, that
\begin{equation} \begin{split}
\label{eq:spatial.error1}
  &\|u^N(t)- u(t)\|_{L^2(\Omega,U)}
  \\
  \leq & \left\|\big(C_N(t) P_N - C(t)\big) u_0\right\|_{L^2(\Omega,U)} + \big\|\big( \Lambda_N^{-\frac{1}{2}}S_N(t)P_N-
  \Lambda^{-\frac{1}{2}}S(t) \big) v_0\big\|_{L^2(\Omega,U)}
  \\
  & + \int_0^t \big\|\Lambda_N^{-\frac{1}{2}}S_N(t-s)P_N F(u^N(s))-\Lambda^{-\frac{1}{2}}S(t-s)F(u(s))  \big\|_{L^2(\Omega,U)} \mbox{d}s
  \\
  & + \Big\| \int_0^t \big(\Lambda_N^{-\frac{1}{2}} S_N(t-s)P_N - \Lambda^{-\frac{1}{2}} S(t-s)\big) \, \mbox{d} W(s)\Big\|_{L^2(\Omega,U)}
  \\
  := & I_1 + I_2 + I_3 + I_4.
\end{split}  \end{equation}
Note that for $\gamma \in [0, \infty)$ and $\xi \in U_N$,
\begin{equation} \label{eq:lambda.comm}
\Lambda_N^{-\gamma} \xi =  \sum_{i=1}^N \lambda_i^{-\gamma} \langle \xi, e_i \rangle e_i = \sum_{i=1}^N \Lambda^{-\gamma} \langle \xi, e_i \rangle e_i = \Lambda^{-\gamma} \xi.
\end{equation}
This together with \eqref{eq:AS.Comm} implies, for all $t \in [0, T]$ and
$\gamma \in [0, \infty)$,  that
%\textcolor{red}{@@}
%\begin{align}
\begin{equation} \begin{split}
\label{eq:Sin.operator2}
  &
  \left\|\Lambda_N^{-\gamma} S_N(t)P_N - \Lambda^{-\gamma} S(t) \right\|_{\mathcal{L}(U)}  = \left\|\Lambda^{-\gamma}S(t) \big(P_N - I\big) \right\|_{\mathcal{L}(U)}
  \\
  =&
  \sup_{i \geq N+1} \big| \lambda_i^{-\gamma} \sin(\sqrt{\lambda_i} t) \big| \leq \lambda_{N+1}^{-\gamma} = \big[\pi^{2} (N+1)^{2}\big]^{-\gamma} \leq N^{-2\gamma}.
\end{split}  \end{equation}
Similarly for all $t \in [0, T]$, $\gamma \in [0, \infty)$, we get
\begin{equation} \begin{split}
\label{eq:Cos.operator1}
  \left\|\Lambda_N^{-\gamma} C_N(t)P_N - \Lambda^{-\gamma} C(t) \right\|_{\mathcal{L}(U)}
  = \sup_{i \geq N+1} \big| \lambda_i^{-\gamma} \cos(\sqrt{\lambda_i} t) \big| \leq \lambda_{N+1}^{-\gamma} \leq N^{-2\gamma}.
\end{split}  \end{equation}
Hence, using \eqref{eq:Sin.operator2} and \eqref{eq:Cos.operator1}  with $\gamma = \frac{\beta}{2} $ shows for all $\beta \in (0, \tfrac{1}{2})$ that
\begin{align}\label{eq:I1I2.estimate}
  I_1 + I_2 \leq \big(\|u_0\|_{L^2(\Omega, \dot{H}^{\beta})} + \|v_0\|_{L^2(\Omega, \dot{H}^{\beta-1})} \big) N^{-\beta}
  \leq  \sqrt{2} \, \|X_0\|_{L^2(\Omega, H^{\beta})} \, N^{-\beta},
\end{align}
where we also used \eqref{eq:AS.Comm} and the fact that $\Lambda^{\gamma}$, $\gamma \in \mathbb{R}$ commutes with $C(t), S(t)$ and $P_N$. With regard to $I_3$, using \eqref{eq:Nf.condition2}, \eqref{eq:AS.Comm}, \eqref{eq:MB.FXFXN}, \eqref{eq:Sin.operator2} and also taking the stability properties of $\Lambda^{-\frac{1}{2}}, S(t)$ and $P_N$ into account, we obtain
\begin{equation} \begin{split}
\label{eq:I3.estimate}
  I_3 \leq& \int_0^t \big\|\Lambda_N^{-\frac{1}{2}}S_N(t-s)P_N \big( F(u^N(s))-F(u(s))\big) \big\|_{L^2(\Omega,U)} \mbox{d}s
  \\
  & + \int_0^t \big\|\big(\Lambda_N^{-\frac{1}{2}}S_N(t-s)P_N -\Lambda^{-\frac{1}{2}}S(t-s)\big)F(u(s))  \big\|_{L^2(\Omega,U)} \mbox{d}s
  \\
  \leq & L \int_0^t \big\|u^N(s)- u(s) \big\|_{L^2(\Omega,U)} \mbox{d} s + \frac{1}{N}\int_0^t \big\|F(u(s))\big\|_{L^2(\Omega, U)}  \mbox{d}s
  \\
  \leq & L \int_0^t \big\|u^N(s)- u(s) \big\|_{L^2(\Omega,U)} \mbox{d} s + \frac{K_3 \, t}{N} \big( \|X_0\|_{L^2(\Omega, H)} + 1 \big).
\end{split}  \end{equation}
Thanks to \eqref{eq:LHS}, \eqref{eq:Lambda.Q}, \eqref{eq:Sin.operator2} and the It\^{o} isometry, one can estimate $I_4$ as follows
\begin{equation} \begin{split}
\label{eq:I4.estimate}
  |I_4|^2 =& \int_0^t \big\|\Lambda_N^{-\frac{1}{2}} S_N(t-s)P_N - \Lambda^{-\frac{1}{2}} S(t-s)\big\|^2_{\mathcal{L}_2(U)} \mbox{d}s
  \\
  \leq & \int_0^t \big\|
  \Lambda^{\frac{\beta-1}{2}} \big\|^2_{\mathcal{L}_2(U)} \cdot \big\|\Lambda^{-\frac{\beta}{2}} S(t-s) \big( P_N - I\big)
  \big\|_{\mathcal{L}(U)}^2
  \mbox{d}s
%  \\
  \leq  \frac{\bar{c}_{\beta} t}{N^{2\beta}}
\end{split}  \end{equation}
for any $\beta < \frac{1}{2}$. Therefore, inserting \eqref{eq:I1I2.estimate}, \eqref{eq:I3.estimate} and \eqref{eq:I4.estimate} into \eqref{eq:spatial.error1} yields
\begin{equation*} \begin{split}
\label{eq:spatial.error2}
  \|u^N(t)- u(t)\|_{L^2(\Omega,U)} \leq&  \sqrt{2} \, \|X_0\|_{L^2(\Omega, H^{\beta})} \, N^{-\beta} + \tfrac{\sqrt{\bar{c}_{\beta} t} } {N^{\beta}} + \tfrac{K_3 t}{N}( \|X_0\|_{L^2(\Omega, H)} + 1 )
  \\
  & + L \smallint_0^t \big\|u^N(s)- u(s) \big\|_{L^2(\Omega,U)} \mbox{d} s
\end{split}  \end{equation*}
for any $\beta < \frac{1}{2}$. Note that $\|u^N(t)- u(t)\|_{L^2(\Omega,U)} \leq
\|u^N(t)\|_{L^2(\Omega,U)} + \|u(t)\|_{L^2(\Omega,U)} < \infty$
by Theorem \ref{thm:unique.mild} and Theorem \ref{thm:Galerkin.unique.mild}. Applying the Gronwall inequality to the preceding estimate with
$\beta = \frac{1}{2}-\epsilon$ gives
\eqref{eq:u.spatial.error}. The proof of Theorem \ref{thm:spatial.error} is thus completed.
$\square$
%Similarly, one can get \eqref{}.

%\section{Further regularity results}
%\label{sect:regularity}

\section{Fully discrete approximations and strong convergence}
\label{sect:superconvergence}

Until now, only spatial discretizations have been investigated. In this section, we turn to the temporal discretizations.
On the interval $[0, T]$, we construct a uniform mesh
$\mathcal{T}_M = \{ t_0, t_1, \cdots, t_M \}$ for $M \in \mathbb{N}$,
satisfying $t_m = m \tau$ with $ \tau = \tfrac{T}{M} $
being the time stepsize.
% More precisely, two exponential time integrators
% (see also \cite{JK09a,JKW11}) will be
% constructed for \eqref{eq:Galerkin.mild.concrete} and
% hence result in full discrete approximations of SWE
% \eqref{eq:SWE}.

\subsection{Fully discrete approximations and main result}

Based on the spatial approximation \eqref{eq:Galerkin.mild.concrete}, we propose two time-stepping schemes as follows:
\begin{align}\label{eq:full.discrete.concrete}
\left\{\!
    \begin{array}{l}
    u^N_{m+1}= C_N(\tau) u_m^N + \Lambda_N^{-\frac{1}{2}} S_N(\tau)v_m^N + \Lambda_N^{-1}\big(I-C_N(\tau) \big)P_N F(u^N_m)
    \\ \:\quad\quad\quad + \int_{t_m}^{t_{m+1}} \Lambda_N^{-\frac{1}{2}} S_N(t_{m+1}-s)P_N \, \mbox{d} W(s),
    \vspace{0.1cm}
    \\
    v_{m+1}^N= -\Lambda_N^{\frac{1}{2}} S_N(\tau) u_m^N + C_N(\tau)v_m^N +  \Lambda_N^{-\frac{1}{2}} S_N(\tau)P_N F(u^N_m)
    \\ \:\quad\quad\quad +  \int_{t_m}^{t_{m+1}} C_N(t_{m+1}-s)P_N \, \mbox{d} W(s),
    \end{array}\right.
\end{align}
and
\begin{align}\label{eq:full.discrete.concrete2}
\left\{\!
    \begin{array}{l}
    u^N_{m+1}= C_N(\tau) u_m^N + \Lambda_N^{-\frac{1}{2}} S_N(\tau)v_m^N + \tau \Lambda_N^{-\frac{1}{2}} S_N(\tau) P_N F(u^N_m)
    \\ \:\quad\quad\quad + \int_{t_m}^{t_{m+1}} \Lambda_N^{-\frac{1}{2}} S_N(t_{m+1}-s)P_N \, \mbox{d} W(s),
    \vspace{0.1cm}
    \\
    v_{m+1}^N= -\Lambda_N^{\frac{1}{2}} S_N(\tau) u_m^N + C_N(\tau)v_m^N +  \tau C_N(\tau)P_N F(u^N_m)
    \\ \:\quad\quad\quad +  \int_{t_m}^{t_{m+1}} C_N(t_{m+1}-s)P_N \, \mbox{d} W(s),
    \end{array}\right.
\end{align}
for $m = 0,1,2, \cdots, M-1$. Here $u^N_{m}$ and $v_{m}^N$ are, respectively, the temporal approximations of $u^N(t)$ and $v^N(t)$ at the grid points $t_m = m \tau$, with the initial values $u^N_0 = P_N u_0$, $v^N_0 = P_N v_0$, and $\tau = T/M$ being the stepsize. It is worthwhile to point out that both proposed schemes are much easier to simulate than it appears at first sight. To show this fact, we take the scheme \eqref{eq:full.discrete.concrete} for example and make some remarks on its implementation.
Observe first that,
%the random variables $\int_{t_m}^{t_{m+1}}\!
%\Lambda_N^{-\frac{1}{2}} S_N(t_{m+1}-s)P_N \,
%\mbox{d} W(s)$ and $\int_{t_m}^{t_{m+1}}  C_N(t_{m+1}-s)
%P_N \, \mbox{d} W(s)$ are common Gaussian distributed
%and independent of $\mathcal{F}_{t_m}$ for all
%$m =0, 1, \cdots, M-1$, $M, N \in \mathbb{N}$.
%As a consequence,
for $i = 1,2, \cdots, N$, $m =0, 1, \cdots, M-1$,
\begin{equation*}
\begin{array}{l}
\zeta_m^i := \left \langle
\int_{t_m}^{t_{m+1}}\!  \Lambda_N^{-\frac{1}{2}} S_N(t_{m+1}-s)P_N \, \mbox{d} W(s), e_i \right\rangle = \lambda_i^{-\frac{1}{2}} \int_{t_m}^{t_{m+1}} \! \sin \! \big((t_{m+1}-s)\lambda_i^{\frac{1}{2}}\big) \, \mbox{d} \beta_i(s)
\end{array}
\end{equation*}
are mutually independent normally distributed random variables satisfying
\begin{align*}
\mathbb{E} \big[ \zeta_m^i \big] =0, \quad  \text{Var}(\zeta_m^i) = \mathbb{E} \big[ |\zeta_m^i|^2 \big] = \tfrac{1}{2\lambda_i} \!\left( \tau-  \tfrac{\sin(2\tau \sqrt{\lambda_i})}{2\sqrt{\lambda_i}} \right).
\end{align*}
Similarly, for $i = 1,2, \cdots, N$, $m =0, 1, \cdots, M-1$,
%\textcolor{red}{@@}
\begin{equation*}
\begin{array}{l}
\widehat{\zeta}_m^i := \left \langle \int_{t_m}^{t_{m+1}}  C_N(t_{m+1}-s) P_N \, \mbox{d} W(s), e_i \right\rangle = \int_{t_m}^{t_{m+1}} \cos \! \big((t_{m+1}-s)\lambda_i^{\frac{1}{2}}\big) \, \mbox{d} \beta_i(s)
\end{array}
\end{equation*}
are mutually independent normally distributed random variables with
\begin{align*}
\mathbb{E} \big[ \widehat{\zeta}_m^i \big] =0, \quad  \text{Var}(\widehat{\zeta}_m^i) = \mathbb{E} \big[ |\widehat{\zeta}_m^i|^2 \big]= \tfrac{1}{2} \! \left( \tau +  \tfrac{\sin(2\tau \sqrt{\lambda_i})}{2\sqrt{\lambda_i}} \right).
\end{align*}
Moreover, the covariance of $\zeta_m^i$ and $\widehat{\zeta}_m^i$ are given by
\begin{equation*}
\text{Cov}(\zeta_m^i, \widehat{\zeta}_m^i) = \mathbb{E} \big[ \zeta_m^i \widehat{\zeta}_m^i \big] = \tfrac{1 - \cos(2 \tau \sqrt{\lambda_i})}{4 \lambda_i} \end{equation*}
for $i = 1,2, \cdots, N$, $m =0, 1, \cdots, M-1$.
Let $D_m^i$  be a family of  $2 \times 2$ matrices with
\begin{equation}\label{eq:covariance.matrix}
Q_m^i  = \left[\! \begin{array}{cc}
     \text{Var}(\zeta_m^i)& \text{Cov}(\zeta_m^i, \widehat{\zeta}_m^i)
    \\
    \text{Cov}(\zeta_m^i, \widehat{\zeta}_m^i) &
    \text{Var}(\widehat{\zeta}_m^i)
    \end{array} \!\right] = D_m^i (D_m^i)^T.
\end{equation}
Then the pair of correlated normally distributed
random variables $(\zeta_m^i, \widehat{\zeta}_m^i)^T$
can be determined by two independent standard normally
distributed random variables:
\begin{equation}
\bigg[ \!\begin{array}{c}
    \zeta_m^i
    \\
    \widehat{\zeta}_m^i
    \end{array}\!\bigg]
    =
    D_m^i
    \bigg[ \!\begin{array}{c}
    R_m^i
    \\
    \widehat{R}_m^i
    \end{array}\!\bigg],
\end{equation}
where $R_m^i \colon \Omega \rightarrow \mathbb{R}$, $\widehat{R}_m^i \colon \Omega \rightarrow \mathbb{R} $ for  $i=1,2,\cdots, N$ and $m =0,1,\cdots, M$ are independent, standard normally distributed random variables.
Accordingly, the components of $u^N_m$ and $v^N_m$ in \eqref{eq:full.discrete.concrete}, i.e., $\langle u^N_m, e_i\rangle$ and $\langle v^N_m, e_i\rangle$ for $i=1,2,\cdots, N$ and $m =0,1,\cdots, M$,  can be calculated by the following recurrence equations:
\begin{align}\label{eq:method.implem.u}
  \left\langle u^N_{m+1}, e_i \right\rangle = & \cos(\tau \sqrt{\lambda_i}) \left\langle u^N_{m}, e_i \right\rangle + \lambda_i^{-\tfrac{1}{2}}\sin(\tau \sqrt{\lambda_i}) \langle v^N_m, e_i\rangle
  \nonumber \\ & +
  \lambda_i^{-1} (1-\cos(\tau \lambda_i^{\tfrac{1}{2}})) \left\langle F(u^N_m), e_i \right\rangle + \zeta_m^i,
  \\
  \left\langle v^N_{m+1}, e_i \right\rangle = & -\lambda_i^{\tfrac{1}{2}} \sin(\tau \sqrt{\lambda_i}) \left\langle u^N_{m}, e_i \right\rangle + \cos(\tau \sqrt{\lambda_i}) \langle v^N_m, e_i\rangle
  \nonumber \\ & +
  \lambda_i^{-\tfrac{1}{2}} \sin(\tau \lambda_i^{\tfrac{1}{2}})  \left\langle F(u^N_m), e_i \right\rangle + \widehat{\zeta}_m^i.
  \label{eq:method.implem.v}
\end{align}
Using the built-in functions "dst" and "idst" in MATLAB,  the scheme \eqref{eq:full.discrete.concrete} can be implemented easily (see Fig. \ref{code} for the implementation code).
Now we formulate our main result as follows.

\begin{theorem}\label{thm:main.result}
Suppose that the nonlinear function $f$ in
\eqref{eq:SWE} satisfies \eqref{f_condition1} and
\eqref{f_condition2}, and let $W(t)$ be the cylindrical $I$-Wiener process represented by
\eqref{W.representation}. Moreover, assume that
%there exists a constant $C \in [0, \infty)$ such that
%\begin{align}\label{eq:init.cond}
$
\|u_0\|_{L^p(\Omega, \dot{H}^1)} +  \|v_0\|_{L^p(\Omega, \dot{H}^0)}  < \infty
$
%\end{align}
for all $p \in [2, 4]$.
Let $u(t)$ be the mild solution of \eqref{eq:SWE} represented by \eqref{eq:mild.concrete} and let $u^N_m$ be the numerical approximation produced by \eqref{eq:full.discrete.concrete} or \eqref{eq:full.discrete.concrete2}, with $\tau = \frac{T}{M}$ being the time stepsize. Then it holds for all $m=0, 1,2,\cdots, M$ and for arbitrarily small $\epsilon >0$ that
\begin{align} \label{eq:main.result}
\left\|u_m^N - u(t_m)\right\|_{L^2(\Omega, U)} \leq K \big(N^{-\frac{1}{2}+\epsilon} + \tau^{1-\epsilon}\big),
\end{align}
where $K \in [0, \infty)$ is a constant depending on $T, \epsilon, L$ and the initial data.
\end{theorem}

The mean-square approximation error
\eqref{eq:main.result} is composed of two parts. The
first term is due to the spatial discretization and the second term is caused by the temporal discretization. The detailed proof of Theorem \ref{thm:main.result}
is postponed to the subsection \ref{sec:proof.main.result}.

\subsection{Some preparatory results}

Before starting the proof of Theorem \ref{thm:main.result}, we need some preparatory results, which are crucial to the convergence analysis.

\begin{lemma} \label{lemma:SCR}
\cite{CLS13,KLL12}
Assume that $S(t)$ and $C(t)$ are the sine and cosine operators as defined above. Then for all $\gamma \in [0,1]$ there exists some constant $\hat{c}\,(\gamma)$ such that
\begin{align}\label{eq:SCNts}
\big\| \big( S(t)- S(s) \big) \Lambda^{-\frac{\gamma}{2}} \big\|_{\mathcal{L}(U)} \leq \hat{c}\, (t-s)^{\gamma},
\;
\big\| \big( C(t)- C(s) \big) \Lambda^{-\frac{\gamma}{2}} \big\|_{\mathcal{L}(U)} \leq \hat{c}\,(t-s)^{\gamma}
\end{align}
for all $t \geq s \geq 0$.
\end{lemma}

Next, a regularity result on the stochastic process $u^N(t)$ is derived, which plays an important role in obtaining the strong convergence rate of the proposed schemes.

\begin{lemma} \label{lemma:Stoch.Conv.Regula.}
Assume that all conditions in Theorem \ref{thm:main.result} are fulfilled.
Then for $\beta \in [0, \frac{1}{2})$ and $\delta \in [0, \tfrac{1}{2}]$ there exists a
constant $K_5(\beta, \delta, p, T) > 0$, depending on $\beta, \delta, p, T$ and
the initial data, such that
%the stochastic convolution $\mathcal{O}^N_t$ given by
%\eqref{eq:stoch.convol.Galerkin} and the stochastic process
%$
%\bar{u}^N(t) :=  u^N(t) - \mathcal{O}^N_t
%$ satisfy
\begin{align}
\left\| u^N(t) - u^N(s) \right\|_{L^p(\Omega, \dot{H}^{-\delta})} \leq& K_5 \, (t-s)^{\beta + \delta},
\quad t > s, \: t, s \in [0, T].
\label{eq:O.U.Holder1}
\end{align}
%for $t > s $ with $t, s \in [0, T]$.
%\textcolor{red}{@ not necessary@ .}
\end{lemma}

{\it Proof.}
Due to \eqref{eq:AS.Comm} and \eqref{eq:lambda.comm}, we first derive from \eqref{eq:Galerkin.mild.concrete} that
\begin{align*}
u^N(t)- u^N(s)
%= &
%C_N(t) P_Nu_0 + \Lambda_N^{-\frac{1}{2}} S_N(t)P_Nv_0
%\nonumber \\ & +
%\smallint_0^t \Lambda_N^{-\frac{1}{2}} S_N(t-r)P_N
%F(u^N(r))\, \mbox{d}r
%+
%\smallint_0^t \Lambda_N^{-\frac{1}{2}} S_N(t-r)P_N \,
%\mbox{d} W(r)
%\nonumber \\ &-
%C_N(s) P_Nu_0 - \Lambda_N^{-\frac{1}{2}} S_N(s)P_Nv_0
%\nonumber \\ & -
%\smallint_0^s \Lambda_N^{-\frac{1}{2}} S_N(s-r)P_N
%F(u^N(r))\, \mbox{d}r - \smallint_0^s
%\Lambda_N^{-\frac{1}{2}} S_N(s-r)P_N \, \mbox{d} W(r)
%\nonumber \\
=&
\big(C(t)-C(s)\big) P_Nu_0 + \Lambda^{-\frac{1}{2}} \big(S(t)-S(s)\big)P_Nv_0
\nonumber \\ & +
\smallint_0^s \Lambda^{-\frac{1}{2}} \big(S(t-r)-S(s-r)\big)P_N F(u^N(r)) \, \mbox{d}r
\nonumber \\& +
\smallint_s^t \Lambda^{-\frac{1}{2}} S(t-r)P_NF(u^N(r)) \, \mbox{d}r
\nonumber \\ & +
\smallint_0^s \Lambda^{-\frac{1}{2}} \big( S(t-r)-S(s-r) \big) P_N \, \mbox{d}W(r)
\nonumber \\ & +
\smallint_s^t \Lambda^{-\frac{1}{2}} S(t-r)P_N \, \mbox{d} W(r).
\end{align*}
Therefore, using the Burkholder-Davis-Gundy type inequality gives
\begin{align*}
& \left\| u^N(t)- u^N(s)\right\|_{L^p(\Omega, \dot{H}^{-\delta})}
\nonumber \\ \leq &
\big\| \Lambda^{-\frac{\delta}{2}}\big(C(t)-C(s)\big) P_Nu_0 \big\|_{L^p(\Omega, U)}
+
\big\|\Lambda^{-\frac{1 + \delta}{2}} \big(S(t)-S(s)\big)P_Nv_0 \big\|_{L^p(\Omega, U)}
\nonumber \\ & +
\smallint_0^s \big\| \Lambda^{-\frac{1 + \delta}{2}} \big(S(t-r)-S(s-r)\big)P_N F(u^N(r)) \big\|_{L^p(\Omega, U)}\mbox{d}r
\nonumber \\& +
\smallint_s^t \big\| \Lambda^{-\frac{1 + \delta}{2}} S(t-r)P_N F(u^N(r))\big\|_{L^p(\Omega, U)} \mbox{d}r
\nonumber \\ & +
\Big( \smallint_0^s \big \|\Lambda^{-\frac{1 + \delta}{2}} \big( S(t-r)-S(s-r) \big) P_N \big\|^2_{\mathcal{L}_2(U)} \, \mbox{d}r \Big)^{1/2}
\nonumber \\ & +
\Big( \smallint_s^t \big\| \Lambda^{-\frac{1 + \delta}{2}} S(t-r)P_N \big\|^2_{\mathcal{L}_2(U)} \, \mbox{d} r \Big)^{1/2}.
\end{align*}
Further, using \eqref{eq:Lambda.Q}, \eqref{eq:MB.FXFXN}, \eqref{eq:SCNts} and the stability properties of $P_N$ and $\Lambda^{-\gamma}, \gamma \geq 0$ results in
\begin{align*}
& \left\|u^N(t)- u^N(s)\right\|_{L^p(\Omega, \dot{H}^{-\delta})}
\nonumber \\ \leq &
\hat{c}\big( \|u_0\|_{L^p(\Omega, \dot{H}^1)} +  \|v_0\|_{L^p(\Omega, U)} \big)(t-s)
\nonumber \\ & +
\hat{c} \smallint_0^s (t-s)\, \|F(u^N(r))\|_{L^p(\Omega, U)} \mbox{d}r
 +
\hat{c} \smallint_s^t (t-r)\, \|F(u^N(r))\|_{L^2(\Omega, U)} \mbox{d}r
\nonumber \\ & +
\Big( \smallint_0^s \big \|\Lambda^{\frac{\beta - 1}{2}} \big\|^2_{\mathcal{L}_2(U)} \big\| \Lambda^{-\frac{\beta + \delta}{2}} \big( S(t-r)-S(s-r) \big) P_N \big\|^2_{\mathcal{L}(U)} \, \mbox{d}r \Big)^{1/2}
\nonumber \\ & +
\Big( \smallint_s^t \big \|\Lambda^{\frac{\beta - 1}{2}} \big\|^2_{\mathcal{L}_2(U)} \big\| \Lambda^{-\frac{\beta + \delta}{2}} S(t-r) P_N \big\|^2_{\mathcal{L}(U)} \, \mbox{d}r \Big)^{1/2}
\nonumber \\ \leq &
\hat{c}\big( \|u_0\|_{L^p(\Omega, \dot{H}^1)} +  \|v_0\|_{L^p(\Omega, U)} \big)(t-s)
\nonumber \\
& +
\hat{c} K_3 ( \|X_0\|_{L^p(\Omega, H)} + 1 ) \big[ T (t-s) + (t-s)^2 \big]
\nonumber \\ & +
\hat{c} \sqrt{ \bar{c}_{\beta} } \big [ \sqrt{T} (t-s)^{\beta + \delta} + (t-s)^{\beta + \delta + \frac{1}{2}} \big]
\nonumber \\ \leq &
K_5 \, (t-s)^{\beta + \delta}
\end{align*}
for $t > s$ with $t, s \in [0, T]$ and where the fact that $\Lambda^{\gamma}$, $\gamma \in \mathbb{R}$ commutes with $C(t), S(t)$ and $P_N$ was also used.  This completes the proof of Lemma \ref{lemma:Stoch.Conv.Regula.}. $\square$

\begin{lemma} \label{lemma:Lambda.F'.Commu.}
Let $F: U \rightarrow U$ be the Nemytskij operator defined by \eqref{eq:Nemytskij}, with the conditions \eqref{f_condition1} and \eqref{f_condition2} fulfilled. Then for $\alpha \in (0, \tfrac{1}{2})$ and $\delta \in (\tfrac{1}{2}, 1]$ there exists a constant $K_6(\alpha, \delta)$ such that
\begin{equation}\label{eq:Lambda.F'.Commu.}
\| F'(\varphi) \psi \|_{\alpha} \leq K_6 \big(\| \varphi \|_{\alpha} + 1\big) \|\psi\|_{\delta}
\end{equation}
holds for all $\varphi \in \dot{H}^{ \alpha}$, $\psi \in \dot{H}^{ \delta} $.
\end{lemma}

{\it Proof. } We denote the space of continuous functions from $(0, 1)$ to $\mathbb{R}$ by $C((0, 1), \mathbb{R})$. The corresponding norm is defined by $\|\varphi\|_{C((0,1), \mathbb{R})} := \sup_{x \in (0,1)} |\varphi(x)|$. Also, we define the Sobolev-Slobodeckij norm (see, e.g., \cite{TV97,TH78}) by
\begin{equation}\label{Sobolev-Slobodeckij-norm}
\|\phi\|_{W^{\alpha, 2}} := \|\phi\| + \left( \int_{0}^1 \int_{0}^1 \frac{|\phi(x) - \phi(y)|^2}{|x-y|^{2 \alpha + 1}} \mbox{d} y \, \mbox{d} x  \right)^{1/2}, \quad \alpha \in (0,1).
\end{equation}
It is known that (see, e.g., (A.46) in \cite{DZ92} or (19.14) in \cite{TV97})
\begin{equation}
\dot{H}^{\alpha} = \{ \phi \in U: \|\phi\|_{ W^{\alpha, 2} }< \infty \}, \quad \mbox{for} \:\: \alpha \in (0, \tfrac{1}{2}),
\end{equation}
and that the norm $\| \cdot \|_{\alpha}$ defined earlier is equivalent to the norm $\|\cdot \|_{W^{\alpha, 2}}$ in $\dot{H}^{\alpha}$:
\begin{equation}\label{eq:norm.equivalent}
\hat{C}_{\alpha} \|\phi\|_{W^{\alpha, 2}} \leq \| \phi \|_{\alpha} \leq C_\alpha \|\phi\|_{W^{\alpha, 2}}, \quad \mbox{for} \:\: \phi \in \dot{H}^{\alpha}, \, \alpha \in (0, \tfrac{1}{2}).
\end{equation}
Then it holds, for all $\varphi \in \dot{H}^{ \alpha}$, $\psi \in \dot{H}^{ \delta} $, $\alpha \in (0, \tfrac{1}{2}), \delta \in (\tfrac{1}{2}, 1]$, that
\begin{align}%\label{eq:ADFA}
& \| F'(\varphi)\psi  \|_{ \alpha } \leq C_\alpha \| F'(\varphi)\psi \|_{W^{\alpha, 2}}
\nonumber \\ =&
C_\alpha \left\|F'(\varphi)\psi \right\|
+ C_\alpha \! \left( \int_0^1\!\!\int_0^1 \frac{\big|\frac{\partial f}{\partial z}(x, \varphi(x))\psi(x)- \frac{\partial f}{\partial z}(y, \varphi(y))\psi(y) \big|^2}{|x-y|^{2\alpha+1 }}\mbox{d}y\mbox{d}x \right)^{\frac{1}{2}}
\nonumber \\ \leq &
LC_\alpha \left\|\psi \right\| + \sqrt{3}C_\alpha \! \left( \int_0^1\!\!\int_0^1 \frac{\big|\frac{\partial f}{\partial z} (x, \varphi(x))\psi(x)- \frac{\partial f}{\partial z}(y, \varphi(x))\psi(x) \big|^2}{|x-y|^{2\alpha+1}} \mbox{d}y\mbox{d}x \right)^{\frac{1}{2}}
\nonumber \\ & +
\sqrt{3}C_\alpha \! \left( \int_0^1\!\!\int_0^1 \frac{\big|\frac{\partial f}{\partial z}(y, \varphi(x))\psi(x)- \frac{\partial f}{\partial z} (y, \varphi(y))\psi(x) \big|^2}{|x-y|^{2\alpha+1}} \mbox{d}y\mbox{d}x \right)^{\frac{1}{2}}
\nonumber \\ & +
\sqrt{3} C_\alpha \! \left( \int_0^1\!\!\int_0^1 \frac{\big|\frac{\partial f}{\partial z} (y, \varphi(y))\big( \psi (x)- \psi(y)\big) \big|^2}{|x-y|^{2\alpha+1}} \mbox{d}y\mbox{d}x \right)^{\frac{1}{2}}
\nonumber \\ \leq &
LC_\alpha \! \left\|\psi \right\| + \sqrt{3}L C_\alpha \! \left( \int_0^1\!\!\int_0^1 |x-y|^{1 - 2 \alpha } \mbox{d}y\mbox{d}x\right)^{\frac{1}{2}} \|\psi\|_{C((0,1), \mathbb{R})}
\nonumber \\ & +
\sqrt{3}L C_\alpha\! \left( \int_0^1\!\!\int_0^1 \frac{\left| \varphi(x)- \varphi(y) \right|^2}{|x-y|^{2\alpha+1}} \mbox{d}y\mbox{d}x \right)^{\frac{1}{2}} \| \psi \|_{C((0,1), \mathbb{R})}
\nonumber \\ & +
\sqrt{3}LC_\alpha \! \left( \int_0^1\!\!\int_0^1 \frac{\left|\psi(x)- \psi(y) \right|^2}{|x-y|^{2\alpha+1}} \mbox{d}y\mbox{d}x \right)^{\frac{1}{2}}
\nonumber \\ \leq &
LC_\alpha\! \left\|\psi \right\| + \sqrt{3}L C_\alpha \| \psi\|_{C((0,1), \mathbb{R})} + \sqrt{3}L C_\alpha \big( \| \varphi \|_{ W^{\alpha,2} }\cdot \| \psi \|_{C((0,1), \mathbb{R})} +  \| \psi \|_{ W^{\alpha,2} }\big)
\nonumber \\ \leq &
LC_\alpha\! \left\|\psi \right\| + \sqrt{3}L C_\alpha \| \psi\|_{C((0,1), \mathbb{R})} + \tfrac{\sqrt{3}L C_\alpha}{ \hat{C}_\alpha} \big( \| \varphi \|_{ \alpha }\cdot \| \psi \|_{ C((0,1), \mathbb{R}) } +  \| \psi \|_{ \alpha }\big)
\nonumber \\ \leq &
K_6 \big(\| \varphi \|_{\alpha} + 1\big) \|\psi\|_{ \delta },
\nonumber
\end{align}
where \eqref{f_condition2}, \eqref{eq:norm.equivalent} and the facts were used that $\dot{H}^{\delta} \subset C((0,1), \mathbb{R})$ continuously for $\delta > \tfrac{1}{2}$ by Sobolev embedding theorem and $\dot{H}^{\alpha} \subset \dot{H}^{\beta}$ for $\alpha \geq \beta$.
$\square$

\subsection{Proof of Theorem \ref{thm:main.result}}
\label{sec:proof.main.result}
To measure the overall mean-square error of the fully discrete schemes,
one can first decompose it as follows:
\begin{align}\label{eq:error.decom}
  \left\|u_m^N - u(t_m)\right\|_{L^2(\Omega, U)} \leq& \left\|u_m^N - u^N(t_m)\right\|_{L^2(\Omega, U)} + \left\|u^N(t_m) - u(t_m)\right\|_{L^2(\Omega, U)}.
\end{align}
Here $m = 0, 1, \cdots, M$ and the second term is the spatial discretization error, which has been estimated in \eqref{eq:u.spatial.error}. As a result, it only remains to estimate the temporal discretization error $\left\|u_m^N - u^N(t_m)\right\|_{L^2(\Omega, U)}$.
We treat the scheme \eqref{eq:full.discrete.concrete} first. By using eigenfunction expansions, one can easily show,  for all $\varphi \in U$, that
\begin{equation}
\begin{array}{ll}
\!\int_{t_m}^{t_{m+1}} \Lambda_N^{-\frac{1}{2}} S_N(t_{m+1}-s) P_N \varphi \, \mbox{d}s = \Lambda_N^{-1} \left(I- C_N(\tau)\right) P_N \varphi, \\
\int_{t_m}^{t_{m+1}} C_N(t_{m+1}-s) P_N \varphi \, \mbox{d}s = \Lambda_N^{-\frac{1}{2}} S_N(\tau) P_N \varphi.
\end{array}
\end{equation}
As a consequence, we can recast \eqref{eq:full.discrete.concrete} as
\begin{align}\label{eq:scheme1.int.form}
\!\! \left\{\!\!
    \begin{array}{l}
    u^N_{m+1}= C_N(\tau) u_m^N + \Lambda_N^{-\frac{1}{2}} S_N(\tau)v_m^N + \int_{t_m}^{t_{m+1}} \! \Lambda_N^{-\frac{1}{2}} S_N(t_{m+1} - s) P_N F(u^N_m) \, \mbox{d} s
    \\ \:\quad\quad\quad + \int_{t_m}^{t_{m+1}} \Lambda_N^{-\frac{1}{2}} S_N(t_{m+1}-s)P_N \, \mbox{d} W(s),
    \vspace{0.1cm}
    \\
    v_{m+1}^N= -\Lambda_N^{\frac{1}{2}} S_N(\tau) u_m^N + C_N(\tau)v_m^N +  \int_{t_m}^{t_{m+1}} C_N(t_{m+1} - s) P_N F(u^N_m) \, \mbox{d} s
    \\ \:\quad\quad\quad +  \int_{t_m}^{t_{m+1}} C_N(t_{m+1}-s)P_N \, \mbox{d} W(s).
    \end{array}\right.
\end{align}
In a compact form, \eqref{eq:scheme1.int.form} is equivalent to
\begin{equation}
\begin{split}
\label{eq:full.dc.compact}
\small
X^N_{m+1} =& E_N(\tau) X_m^N + \!\int_{t_m}^{t_{m+1}}\! E_N(t_{m+1}-s)\,\mathbf{F}_N(X_m^N) \,\mbox{d}s
\\ &
+ \!\int_{t_m}^{t_{m+1}}\! E_N(t_{m+1}-s) B_N \mbox{d} W(s),
\end{split}  \end{equation}
where we denote
$
X_m^N = (
    u_m^N, \,
    v_m^N)^T.
$
for $m = 0, 1, \cdots, M$.
This further implies that
\begin{equation} \begin{split}
\label{eq:full.dc.rec}
X^N_{m+1} =& E_N(t_{m+1}) X_0^N + \sum_{l=0}^{m} \int_{t_l}^{t_{l+1}} E_N(t_{m+1}-s) \mathbf{F}_N(X^N_l) \, \mbox{d} s
\\
&+ \int_0^{t_{m+1}} E_N(t_{m+1}-s) B_N \, \mbox{d} W(s).
\end{split}  \end{equation}
Subtracting \eqref{eq:Galerk.mild.solut} from \eqref{eq:full.dc.rec} yields
\begin{equation}\label{eq:X_Diff}
%\begin{array}{l}
X_{m+1}^N - X^N(t_{m+1}) = \sum\limits_{l=0}^m \int_{t_l}^{t_{l+1}} E_N(t_{m+1}-s) \Big(\mathbf{F}_N(X_l^N)-\mathbf{F}_N(X^N(s))\Big) \mbox{d}s,
%\end{array}
\end{equation}
which in turn implies, for all $m = 0,1,\cdots, M-1$, that
\begin{equation}\label{eq:u.time.diff}
\begin{split}
u_{m+1}^N - u^N(t_{m+1}) = & \sum_{l=0}^m \!\int_{t_l}^{t_{l+1}} \Lambda_N^{-\frac{1}{2}} S_N(t_{m+1}-s) P_N \! \left( F(u^N_l)-F(u^N(s))\right) \mbox{d} s
\\ = &
\sum_{l=0}^m \!\int_{t_l}^{t_{l+1}} \Lambda^{-\frac{1}{2}} S(t_{m+1}-s) P_N \!\left( F(u^N_l)-F(u^N(s))\right) \mbox{d} s.
\end{split}
\end{equation}
%and
%\begin{align}\label{eq:v.time.diff}
%v_{m+1}^N - v^N(t_{m+1}) = \sum_{l=0}^m \int_{t_l}^{t_{l+1}}
%C_N(t_{m+1}-s) P_N \left( F(u^N_l)-F(u^N(s))\right)
%\mbox{d} s
%\end{align}
Using \eqref{eq:Lambda.Galerkin}, \eqref{eq:AS.Comm}, \eqref{eq:Nf.condition2}, the stabilities of $\Lambda^{-\frac{1}{2}}, P_N$ and $S(t)$ shows that
\begin{align}\label{eq:time.error.estimate}
&\|u^N_{m+1}-u^N(t_{m+1}) \|_{L^2(\Omega, U)}
%\nonumber \\ = &
%\bigg\|\sum_{l=0}^m \int_{t_l}^{t_{l+1}}
%\Lambda^{-\frac{1}{2}} S(t_{m+1}-s) P_N \left( F(u^N_l)
%-F(u^N(s))\right) \mbox{d} s\bigg\|_{L^2(\Omega, U)}
\nonumber \\ \leq &
\sum_{l=0}^m \int_{t_l}^{t_{l+1}}
\big\|\Lambda^{-\frac{1}{2}} S(t_{m+1}-s) P_N
\left( F(u^N_l)-F(u^N(t_l))\right)
\big\|_{L^2(\Omega, U)} \mbox{d} s
%\nonumber \\ &
+ J
%\bigg\|\sum_{l=0}^m \int_{t_l}^{t_{l+1}}
%\Lambda^{-\frac{1}{2}} S(t_{m+1}-s) P_N \left( F(u^N(s))
%- F(u^N(t_l))\right) \mbox{d} s\bigg\|_{L^2(\Omega, U)}
\nonumber \\ \leq &
\sum_{l=0}^m \int_{t_l}^{t_{l+1}} \left\| F(u^N_l)-F(u^N(t_l)) \right\|_{L^2(\Omega, U)} \mbox{d} s + J
\nonumber \\ \leq &
L \tau \sum_{l=0}^m  \left\| u^N_l- u^N(t_l) \right\|_{L^2(\Omega, U)}  + J,
\end{align}
where for simplicity we denote
\begin{equation}\label{eq:estimate.J}
\begin{array}{l}
  J = \Big\|\sum\limits_{l=0}^m \int_{t_l}^{t_{l+1}} \Lambda^{-\frac{1}{2}} S(t_{m+1}-s) P_N \big( F(u^N(s)) - F(u^N(t_l))\big) \mbox{d} s\Big\|_{L^2(\Omega, U)}.
\end{array}
\end{equation}
Now it remains to estimate $J$. Thanks to the stability properties of $S(t), P_N$ and the fact that $\Lambda^{-1/2}$ commutes with $S(t), P_N$, we get
\begin{equation} \begin{split}
\label{eq:J.first.estimate}
J \leq &
\sum_{l=0}^m \int_{t_l}^{t_{l+1}}
\big\| \Lambda^{-\frac{1}{2}} S(t_{m+1}-s) P_N \big( F(u^N(s))-F(u^N(t_l)) \big) \big\|_{L^2(\Omega, U)} \mbox{d} s
\\
\leq &
\sum_{l=0}^m \int_{t_l}^{t_{l+1}}
\big\| \Lambda^{-\frac{1}{2}} \big( F(u^N(s))-F(u^N(t_l)) \big) \big\|_{L^2(\Omega, U)} \mbox{d} s
\\
= &
\sum_{l=0}^m \int_{t_l}^{t_{l+1}}
\Big\| \!\int_0^1 \Lambda^{-\frac{1}{2}} F'\big(  \chi(t_l,s, r) \big) \big( u^N(s) -u^N(t_l) \big) \mbox{d} r \Big\|_{L^2(\Omega, U)} \mbox{d} s
\\
\leq &
\sum_{l=0}^m \int_{t_l}^{t_{l+1}}\!
\int_0^1 \big\| \Lambda^{-\frac{1}{2}} F'\big( \chi(t_l,s, r)  \big) \big( u^N(s) -u^N(t_l) \big) \big\|_{L^2(\Omega, U)} \mbox{d} r \, \mbox{d} s,
\end{split}
\end{equation}
where for simplicity of presentation we denote
\begin{equation}
\chi(t_l,s, r) := u^N(t_l) + r(u^N(s) -u^N(t_l)).
\end{equation}
In the next step we start to estimate the last term of \eqref{eq:J.first.estimate}:
\begin{equation}\label{eq:J22}
\begin{split}
& \big\| \Lambda^{-\frac{1}{2}} F'\big( \chi(t_l,s, r)  \big) \big( u^N(s) -u^N(t_l) \big) \big\|_{L^2(\Omega, U)}
\\ = &
\Big\| \big \|\Lambda^{-\frac{1}{2}} F'\big( \chi(t_l,s, r)  \big) \big( u^N(s) -u^N(t_l) \big) \big\| \Big\|_{L^2(\Omega, \mathbb{R})}
\\ = &
\Big\|\sup_{ \psi \in U, \|\psi\|\leq 1 } \big|\big\langle \Lambda^{-\frac{1}{2}} F'\big( \chi(t_l,s, r)\big) \big( u^N(s) -u^N(t_l) \big), \, \psi\big\rangle \big|\,\Big\|_{L^2(\Omega, \mathbb{R})}
\\ = &
\Big\|\sup_{ \psi \in U, \|\psi\|\leq 1} \big|\big\langle u^N(s) -u^N(t_l), \, \big(F'\big( \chi(t_l,s, r)\big)\big)^* \Lambda^{-\frac{1}{2}}  \psi\big\rangle \big|\,\Big\|_{L^2(\Omega, \mathbb{R})}
\\ = &
\Big\|\sup_{\psi \in U, \|\psi\|\leq 1} \big|\big\langle \Lambda^{\frac{\epsilon - 1 }{4}} \big( u^N(s) -u^N(t_l) \big), \,
\Lambda^{\frac{1-\epsilon}{4}} F'\big( \chi(t_l,s, r)\big)  \Lambda^{-\frac{1}{2}}\psi\big\rangle \big|\,\Big\|_{L^2(\Omega, \mathbb{R})}
\\
\leq &
\Big\| \big\| \Lambda^{ \frac{\epsilon- 1}{4}}
\big( u^N(s) -u^N(t_l) \big) \big\|
\sup_{\psi \in U, \|\psi\|\leq 1} \big\|\Lambda^{\frac{1-\epsilon}{4}}  F'( \chi(t_l,s, r) )\Lambda^{-\frac{1}{2}}\psi \big\|\,\Big\|_{L^2(\Omega, \mathbb{R})}
\\
\leq &
\big\|
   u^N(s) -u^N(t_l)
  \big\|_{L^4(\Omega, \dot{H}^{\frac{\epsilon-1}{2}})} \,
\Big\|\sup_{\psi \in U, \|\psi\|\leq 1} \! \big\| F'( \chi(t_l,s, r) )\Lambda^{-\frac{1}{2}}\psi \big\|_{\frac{1-\epsilon}{2}} \Big\|_{L^4(\Omega, \mathbb{R})},
\end{split}
\end{equation}
where $\epsilon \in (0, 1)$ is arbitrarily small and where the stability properties of $\Lambda^{-\frac{1}{2}},S(t),P_N$, self-adjointness of $F'(u)$ and $\Lambda^{\gamma}, \gamma \in \mathbb{R}$,  the Cauchy-Schwarz inequality and H\"{o}lder's inequality were invoked. In view of \eqref{eq:XN.MB} and \eqref{eq:Lambda.F'.Commu.} with $\alpha = \tfrac{1-\epsilon}{2}, \delta = 1$ we get
\begin{align}
& \Big\| \sup_{\psi \in U, \|\psi\|\leq 1} \big\|  F'\big( \chi(t_l,s, r)  \big)\Lambda^{-\frac{1}{2}}\psi \big\|_{\frac{1-\epsilon}{2}}
\Big\|_{L^4(\Omega, \mathbb{R})} \leq  K_6 \left\| \big\| \chi(t_l,s, r) \big\|_{\frac{1-\epsilon}{2}} + 1 \right\|_{L^4(\Omega, \mathbb{R})}
\nonumber \\ \leq &
K_6 \Big( \sup_{s \in [0, T]} \|u^N(s)\|_{L^4(\Omega, \dot{H}^{\frac{1-\epsilon}{2}})} + 1 \Big)
\leq
K_6  \big( K_2 \|X_0\|_{L^4(\Omega, H^{\frac{1-\epsilon}{2}})} + K_2 + 1 \big).
\label{eq:Last.estimate.J22}
\end{align}
Furthermore, using \eqref{eq:O.U.Holder1} with $\beta = \delta = \frac{1 - \epsilon}{2} $ gives
\begin{equation}\label{eq:time.regul.negative.Sobelev}
\big\|
  u^N(s) -u^N(t_l)
  \big\|_{L^4(\Omega, \dot{H}^{\frac{\epsilon-1}{2}})} \leq K_5 (s - t_l)^{1 - \epsilon}, \quad s \in [t_l, t_{l+1}].
\end{equation}
Thus plugging the above two estimates into \eqref{eq:J22} shows, for arbitrarily small $\epsilon \in (0,1)$
\begin{align}\label{eq:J22.estimate.final}
  \big\| \Lambda^{-\frac{1}{2}} F'\big( \chi(t_l,s, r)  \big) \big( u^N(s) -u^N(t_l) \big) \big\|_{L^2(\Omega, U)}  \leq  K_7 \, \tau^{1-\epsilon},
\end{align}
where $K_7 =  K_5 K_6 ( K_2\|X_0\|_{L^4(\Omega, H^{(1-\epsilon)/2})} + K_2 + 1 )$.
This together with \eqref{eq:J.first.estimate} implies
\begin{align}\label{eq:J.final}
J \leq  T K_7  \, \tau^{1-\epsilon}
\end{align}
for arbitrarily small $\epsilon \in (0,1)$.
Hence one can deduce from \eqref{eq:time.error.estimate} that, for arbitrarily small $\epsilon \in (0,1)$ and all $m = 0, 1, \cdots, M- 1$
\begin{align}\label{eq:time.error.estimate2}
\|u^N_{m+1}-u^N(t_{m+1}) \|_{L^2(\Omega, U)} \leq L \tau \sum_{l=0}^m  \left\| u^N_l- u^N(t_l) \right\|_{L^2(\Omega, U)}  + T K_7  \, \tau^{1-\epsilon}.
\end{align}
%This recurrence equation obviously implies
%$\|u^N_{m}-u^N(t_{m}) \|_{L^2(\Omega, U)} < \infty$
%for all $m =1,2,\cdots, M$ and thus
The discrete Gronwall inequality applied to \eqref{eq:time.error.estimate2} gives, for all $m =1,2,\cdots, M$ and arbitrarily small $\epsilon \in (0,1)$, that
\begin{align}\label{eq:time.error.estimate.final}
\left\|u_m^N - u^N(t_m)\right\|_{L^2(\Omega, U)}  \leq T K_7 e^{LT} \,\tau^{1-\epsilon}.
\end{align}
Finally, this together with \eqref{eq:u.spatial.error}, \eqref{eq:error.decom} and the assumption on initial data gives the overall discretization
error \eqref{eq:main.result} for the first numerical scheme \eqref{eq:full.discrete.concrete}.

For the numerical scheme \eqref{eq:full.discrete.concrete2}, one can also write it in a compact form:
\begin{equation}
\label{eq:full.dc.compact2}
\begin{array}{l}
X^N_{m+1} = E_N(\tau) X_m^N + \tau E_N(\tau) \mathbf{F}_N(X_m^N)
+ \int_{t_m}^{t_{m+1}}\! E_N(t_{m+1}-s) B_N \mbox{d} W(s).
\end{array}
\end{equation}
Similarly to \eqref{eq:X_Diff}, one has
\begin{equation*}
X_{m+1}^N - X^N(t_{m+1}) = \sum_{l=0}^m \smallint_{t_l}^{t_{l+1}}  \left( E_N(t_{m+1}-t_l) \mathbf{F}_N(X_l^N)- E_N(t_{m+1}-s) \mathbf{F}_N(X^N(s))\right) \mbox{d}s.
\end{equation*}
This shows, for all $m = 0,1,2, \cdots, M-1$, that
\begin{align}
& u_{m+1}^N - u^N(t_{m+1})
\nonumber \\ = &
\sum_{l=0}^m \int_{t_l}^{t_{l+1}} \!\left( \Lambda^{-\frac{1}{2}} S(t_{m+1}-t_l) P_N F(u^N_l)- \Lambda^{-\frac{1}{2}} S(t_{m+1}-s) P_N F(u^N(s))\right) \mbox{d} s
\nonumber \\ =&
\sum_{l=0}^m \int_{t_l}^{t_{l+1}} \!
\Lambda^{-\frac{1}{2}} S(t_{m+1}-t_l) P_N \Big(F(u^N_l)- F(u^N(s))\Big)\, \mbox{d} s
\nonumber \\& +
\sum_{l=0}^m \int_{t_l}^{t_{l+1}} \!\left( \Lambda^{-\frac{1}{2}} S(t_{m+1}-t_l) - \Lambda^{-\frac{1}{2}} S(t_{m+1}-s) \right) P_N F(u^N(s)) \,\mbox{d} s.
\label{eq:u.time.diff2}
\end{align}
Therefore, it follows that, for all $m = 0,1,2, \cdots, M-1$
\begin{align}
& \|u_{m+1}^N - u^N(t_{m+1})\|_{L^2(\Omega, U)}
\nonumber \\ \leq&
\bigg\|\sum_{l=0}^m \int_{t_l}^{t_{l+1}} \Lambda^{-\frac{1}{2}} S( t_{m+1}-t_l ) P_N \Big(F(u^N_l)- F(u^N(s))\Big) \mbox{d}s \bigg\|_{L^2(\Omega, U)}
\nonumber \\& +
\sum_{l=0}^m \int_{t_l}^{t_{l+1}} \!\left\|\left( \Lambda^{-\frac{1}{2}} S( t_{m+1}-t_l ) - \Lambda^{-\frac{1}{2}} S(t_{m+1}-s) \right) P_N F(u^N(s))\right\|_{L^2(\Omega, U)} \mbox{d} s
\nonumber \\ \leq&
\bigg\|\sum_{l=0}^m \int_{t_l}^{t_{l+1}}\! \Lambda^{-\frac{1}{2}} S( t_{m+1}-t_l ) P_N \Big(F(u^N_l)- F(u^N(t_l))\Big) \mbox{d}s \bigg\|_{L^2(\Omega, U)} + J'
\nonumber \\ & +
\sum_{l=0}^m \int_{t_l}^{t_{l+1}} \hat{c} (s - t_l) \, \|F(u^N(s))\|_{L^2(\Omega, U)} \mbox{d} s
\nonumber \\ \leq&
\sum_{l=0}^m \int_{t_l}^{t_{l+1}} \Big\| \Lambda^{-\frac{1}{2}} S( t_{m+1}-t_l ) P_N \Big(F(u^N_l)- F(u^N(t_l))\Big) \Big\|_{L^2(\Omega, U)} \mbox{d}s + J'
\nonumber \\ & +
\hat{c} K_3 T ( \|X_0\|_{L^p(\Omega, H)} + 1 ) \tau
\nonumber \\ \leq&
L \tau \sum_{l=0}^m  \big\|u^N_l - u^N(t_l)\big \|_{L^2(\Omega, U)} + J' +  \hat{c} K_3 T ( \|X_0\|_{L^p(\Omega, H)} + 1 ) \tau,
\label{eq:scheme2.error}
\end{align}
where we denote
\begin{equation}
\begin{array}{l}
J' = \Big\|\sum\limits_{l=0}^m \int_{t_l}^{t_{l+1}} \Lambda^{-\frac{1}{2}} S( t_{m+1}-t_l ) P_N \! \left( F(u^N(s)) - F(u^N(t_l))\right) \mbox{d} s\Big\|_{L^2(\Omega, U)}.
\end{array}
\end{equation}
Here \eqref{eq:Nf.condition2}, \eqref{eq:MB.FXFXN}, \eqref{eq:SCNts} and the stability properties of $\Lambda^{-\frac{1}{2}},S(t),P_N$ were also used in these estimates.
Note that the term $J'$ is almost a copy of the term $J$ in \eqref{eq:estimate.J},
with only $S_N(t_{m+1}-s)$ replaced by $S_N( t_{m+1}-t_l )$.
Thus, repeating exactly the same lines as before one can easily obtain for arbitrarily small $\epsilon \in (0,1)$ that
\begin{align}\label{eq:J'.final}
J' \leq T K_7 \, \tau^{1-\epsilon}.
\end{align}
Inserting \eqref{eq:J'.final} into \eqref{eq:scheme2.error}, we apply the discrete Gronwall inequality to obtain
the temporal discretization error for the scheme \eqref{eq:full.discrete.concrete2}.
Finally, taking  \eqref{eq:u.spatial.error} and \eqref{eq:error.decom} into account finishes the proof of Theorem \ref{thm:main.result}.
$\square$

%\section{Extensions ...}
%\label{sect:extensions}

\section{Numerical results}
\label{sect:implementation}

In this section, we perform some numerical
examples to illustrate our previous findings. As the first numerical example, we consider the Sine-Gordon equation driven by space-time white noise:
\begin{equation}\label{eq:SWE.nonlinear1}
\left\{
    \begin{array}{lll}
    \frac{\partial^2 u}{\partial t^2} = \frac{\partial^2 u}{\partial x^2} -\sin( u ) + \dot{W}, \quad  t \in (0, 1], \:\: x \in (0,1),
    \\
     u(0, x) = \frac{\partial u}{\partial t}  (0,x) = 0, \: x \in (0,1),
     \\
     u(t, 0) = u(t,1) = 0,  \: t >0.
    \end{array}\right.
\end{equation}
The corresponding deterministic equation is used to
describe the dynamics of coupled Josephson junctions
driven by a fluctuating current source \cite{Cpl07}.
Next, we use various numerical schemes to
solve \eqref{eq:SWE.nonlinear1} and compare their
computational errors. Note that the expectations are approximated by computing averages over $100$ samples.

The MATLAB code of one-path simulation of \eqref{eq:SWE.nonlinear1} using the scheme \eqref{eq:full.discrete.concrete} is presented in
Fig.~\ref{code}. Here we evoke the built-in functions "dst" and "idst" in MATLAB, which are based on
the fast Fourier transform, to numerically approximate the inner products $\langle F(u^N_m), e_i \rangle$ in \eqref{eq:method.implem.u}-\eqref{eq:method.implem.v} at cheap costs. The aliasing errors are then neglected. More details and remarks on such implementation can be found in \cite[Section 5.1]{WG13}.

\begin{figure}[htp]
\centering
\begin{verbatim}
M = 64; N = M^2; T = 1; tau = T/M;
A = pi^2*(1:N).^2;  sqrtA = sqrt(A);
CosA = cos(sqrtA*tau); SinA = sin(sqrtA*tau);
Var1 = (tau-sin(2*sqrtA*tau)./(2*sqrtA))./(2*A);
Var2 = (tau+sin(2*sqrtA*tau)./(2*sqrtA))./2;
Cov12 = sin(sqrtA*tau).^2./(2*A);
SW1 = sqrt(Var1); SW21 = Cov12./sqrt(Var1);
SW22 = sqrt( (Var1.*Var2 - Cov12.^2)./Var1 );
f = @(x) -sin(x);
Y1 = zeros(1,N); Y2 = zeros(1,N);
for m = 1:M
    Rd1 = randn(1,N); Rd2 = randn(1,N);
    dW1 = SW1.*Rd1; dW2 = SW21.*Rd1 + SW22.*Rd2;
    y1 = dst(Y1)*sqrt(2);
    Fy1 = idst(f(y1))/sqrt(2);
    Y10 = Y1;
    Y1 = CosA.*Y1 + 1./sqrtA.*SinA.*Y2 + 1./A.*(1-CosA).*Fy1 + dW1;
    Y2 = -sqrtA.*SinA.*Y10 + CosA.*Y2 + 1./sqrtA.*SinA.*Fy1 + dW2;
end
plot((0:N+1)/(N+1),[0,dst(Y1)*sqrt(2),0]);
\end{verbatim}
\caption{MATLAB code for one-path simulation of \eqref{eq:SWE.nonlinear1} using \eqref{eq:full.discrete.concrete}}
\label{code}
\end{figure}

\begin{figure}[htp]
         \centering
         \includegraphics[width=3in,height=2.5in]
         {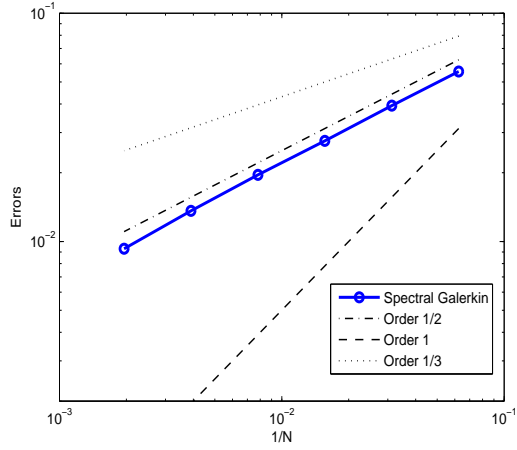}
         \caption{Spatial errors for the spectral Galerkin method applied to SWE \eqref{eq:SWE.nonlinear1}.}
         \label{fig:spatial.error}
\end{figure}

Now let us start with tests on the convergence rates. At first, we consider the spatial convergence rate of the spectral Galerkin method \eqref{eq:Galerkin.SEE}. Fig. \ref{fig:spatial.error} depicts the spatial approximation errors $\|u(T)-u^N(T)\|_{L^2(\Omega, U)}$ against $\frac{1}{N}$ on a log-log scale, with $T = 1$ and $N = 2^j, j =4,5,...,9$. One can detect that the errors decrease at a slope of order $\frac{1}{2}$ as $\frac{1}{N}$ decreases, which is consistent with our previous assertion on the spatial convergence rate. Note that for the temporal discretization we used here the proposed scheme \eqref{eq:full.discrete.concrete2} at a small time stepsize $\tau_{exact} = 2^{-7}$. In addition, $N_{exact} = 2^{14}$ was used to compute the ``exact'' solution $u(T)$.

%To this end,  "exact" solutions
%such as $u(t)$ and $u^N(t)$ are always needed,
%which in the following numerical results are identified
%with the considered numerical solutions
%using very small stepsize.

%(here and below we write $b-$ for the convergence order
%if the convergence order is higher than $b-\epsilon$
%for arbitrarily small $0< \epsilon <b$).

\begin{figure}[htp]
         \centering
         \includegraphics[width=2.5in,height=2.3in]
         {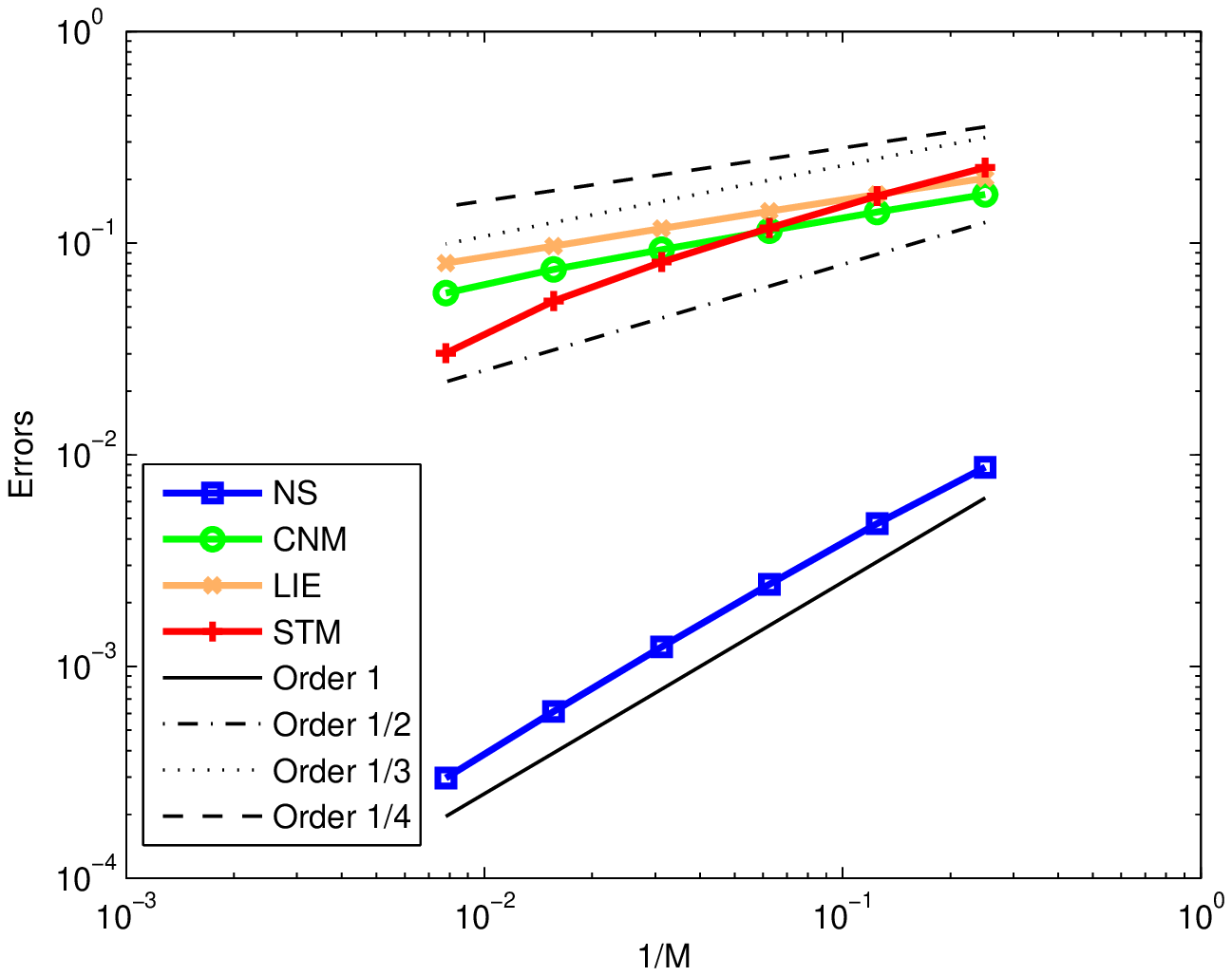}
         \includegraphics[width=2.5in,height=2.3in]
         {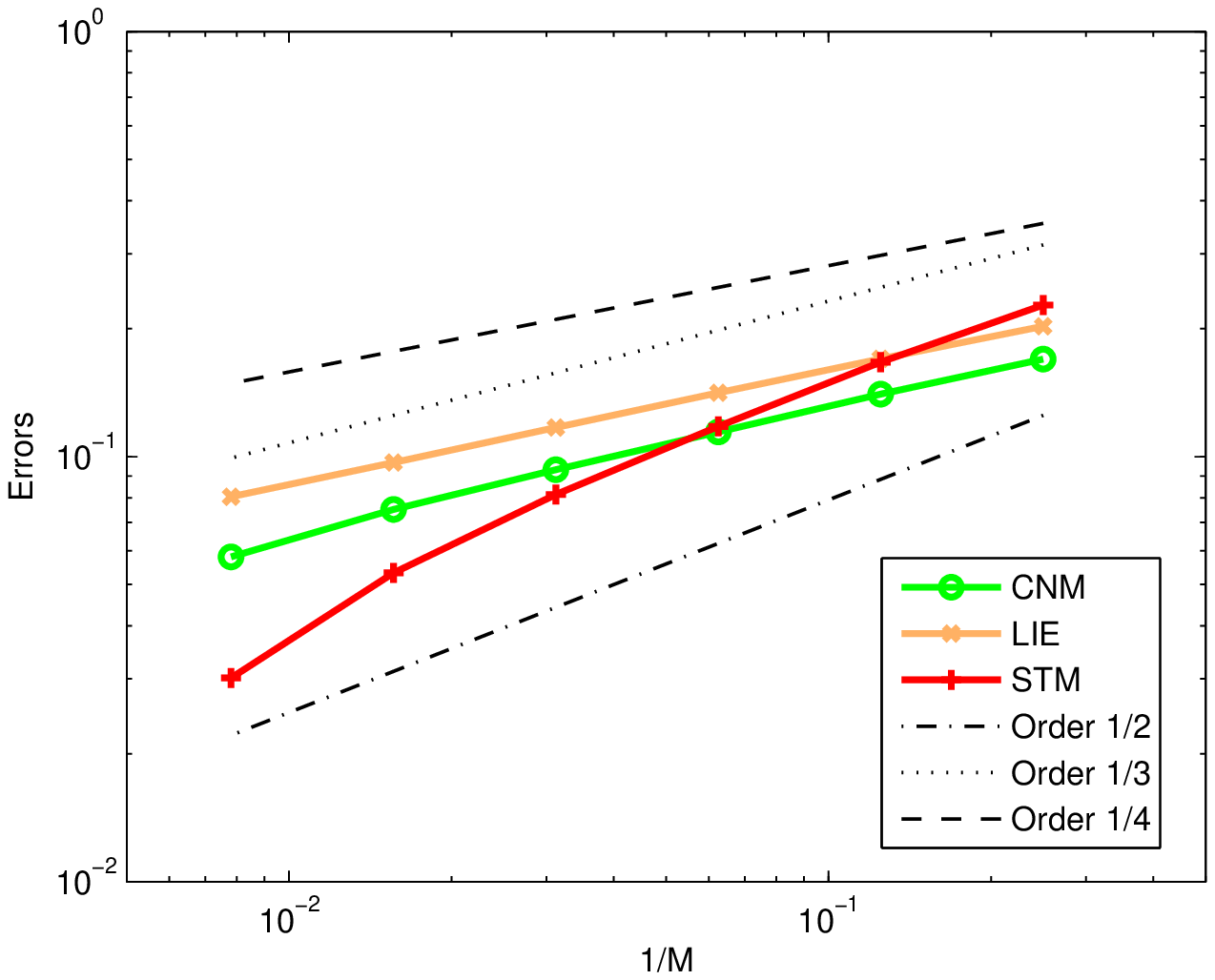}
         \caption{Strong convergence rates in time for various time integrators applied to \eqref{eq:Galerkin.SEE}.}
         \label{fig:time.error}
\end{figure}

%$M = 2^8; N = 2^{12}; rounds = 100; T = 1;
%tau = T/M; Ns= 2^{(4:9)}; $
We now fix $N = 100$ and compare several time integrators applied to approximate $u^N(T),\, T=1$ using various stepsizes $\tau = 2^{-j}, j =2,3,...,7$. Again, the ``exact'' solution $u^N(1)$ is approximated by the method \eqref{eq:full.discrete.concrete2} with a very small stepsize $\tau_{exact} = 2^{-10}$. In Fig.~\ref{fig:time.error}, we present approximation errors caused by different temporal discretizations, including the linear implicit Euler (LIE) scheme \cite{KLL12}, the Crank-Nicolson-Maruyama (CNM) scheme \cite{HE03,KLL12}, the stochastic trigonometric method (STM)\cite{CLS13,Wang13} and the proposed scheme \eqref{eq:full.discrete.concrete}. From the left picture of Fig.~\ref{fig:time.error}, one can easily observe that the scheme \eqref{eq:full.discrete.concrete} performs much better than the other ones. On the one hand, it produces much smaller errors than the other three schemes. On the other hand, computational errors of the scheme \eqref{eq:full.discrete.concrete} decrease faster, i.e., with rate $1$. To clearly display the convergence rates of the three existing schemes, we change the scales of coordinate axes and hide the computational errors of the scheme \eqref{eq:full.discrete.concrete}. In the right picture of Fig.~\ref{fig:time.error}, one can find that the approximation errors of the linear implicit Euler (LIE) scheme decrease with order $\frac{1}{4}$, the errors of the Crank-Nicolson-Maruyama (CNM) scheme with order $\frac{1}{3}$, and the errors of the stochastic trigonometric method (STM) with order $\frac{1}{2}$. These numerical performances coincide with theoretical findings, see Theorem \ref{thm:main.result}, \cite[Theorem 4.12]{KLL12} and \cite[Theorem 4.1]{CLS13}.

%$N = 100; M = 2^{11}; Ms = 2.^{(2:7)}; rounds = 200;$

\begin{figure}[htp]
         \centering
         \includegraphics[width=2.5in,height=2.3in]
         {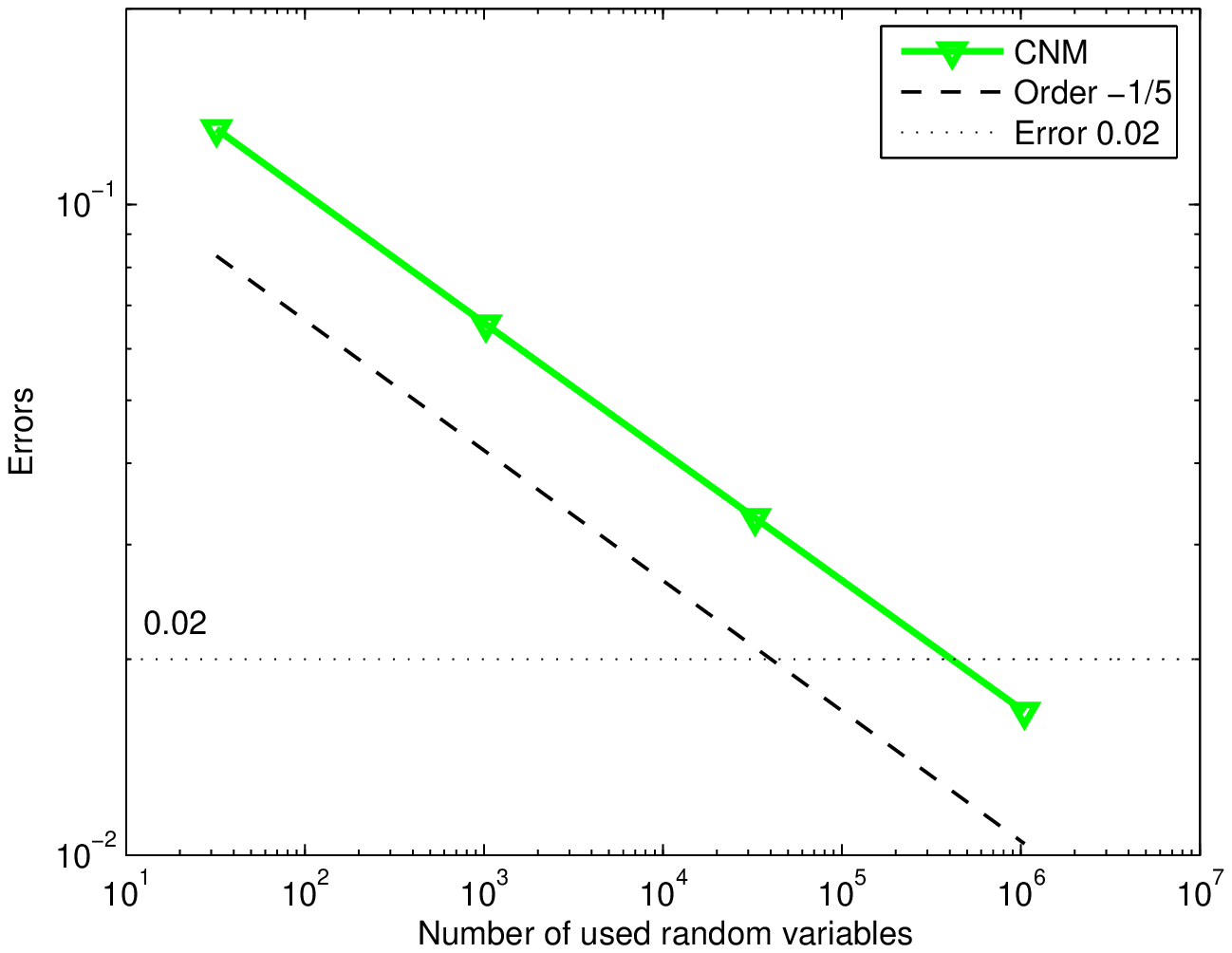}
         \includegraphics[width=2.5in,height=2.3in]
         {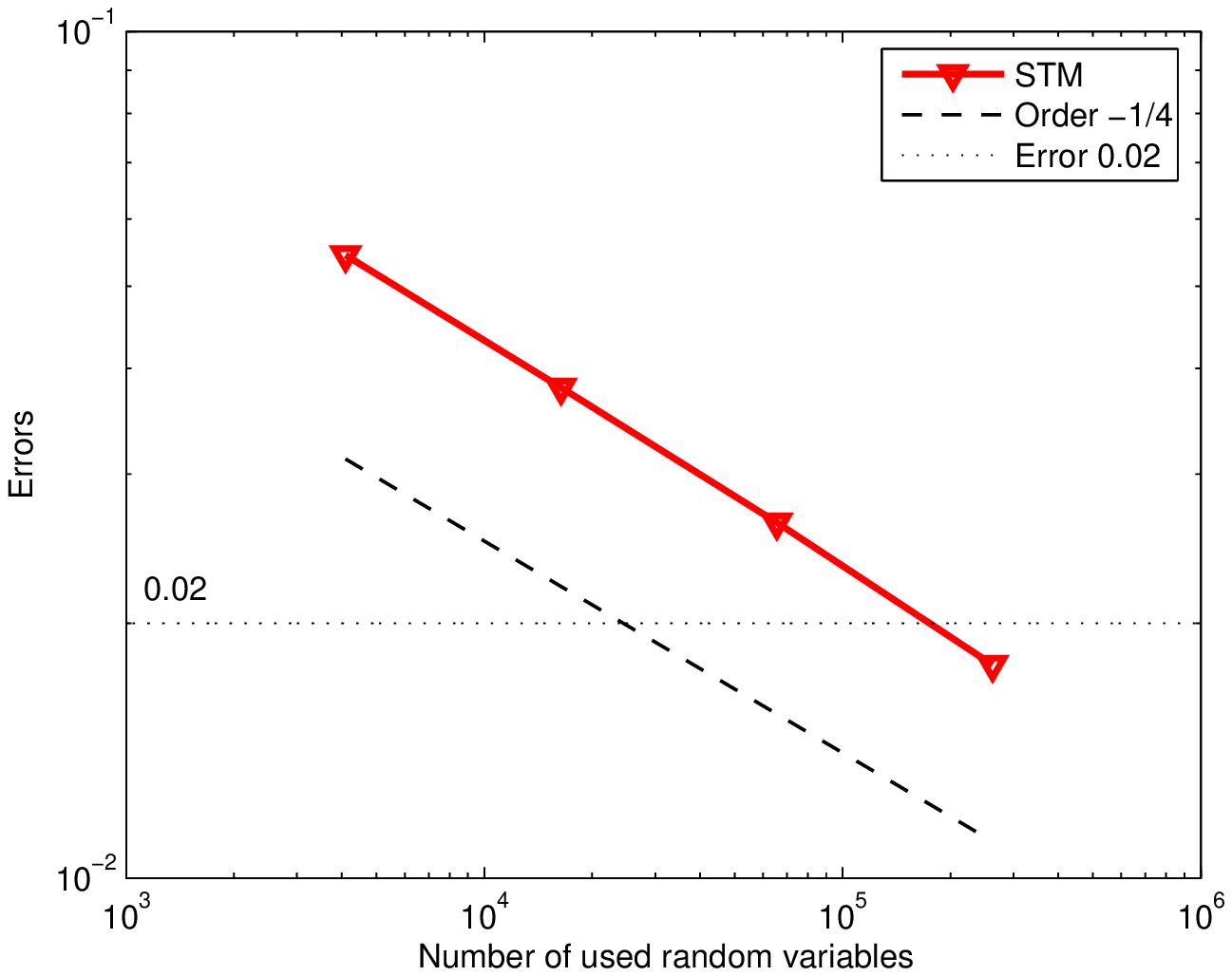}
         \includegraphics[width=2.5in,height=2.3in]
         {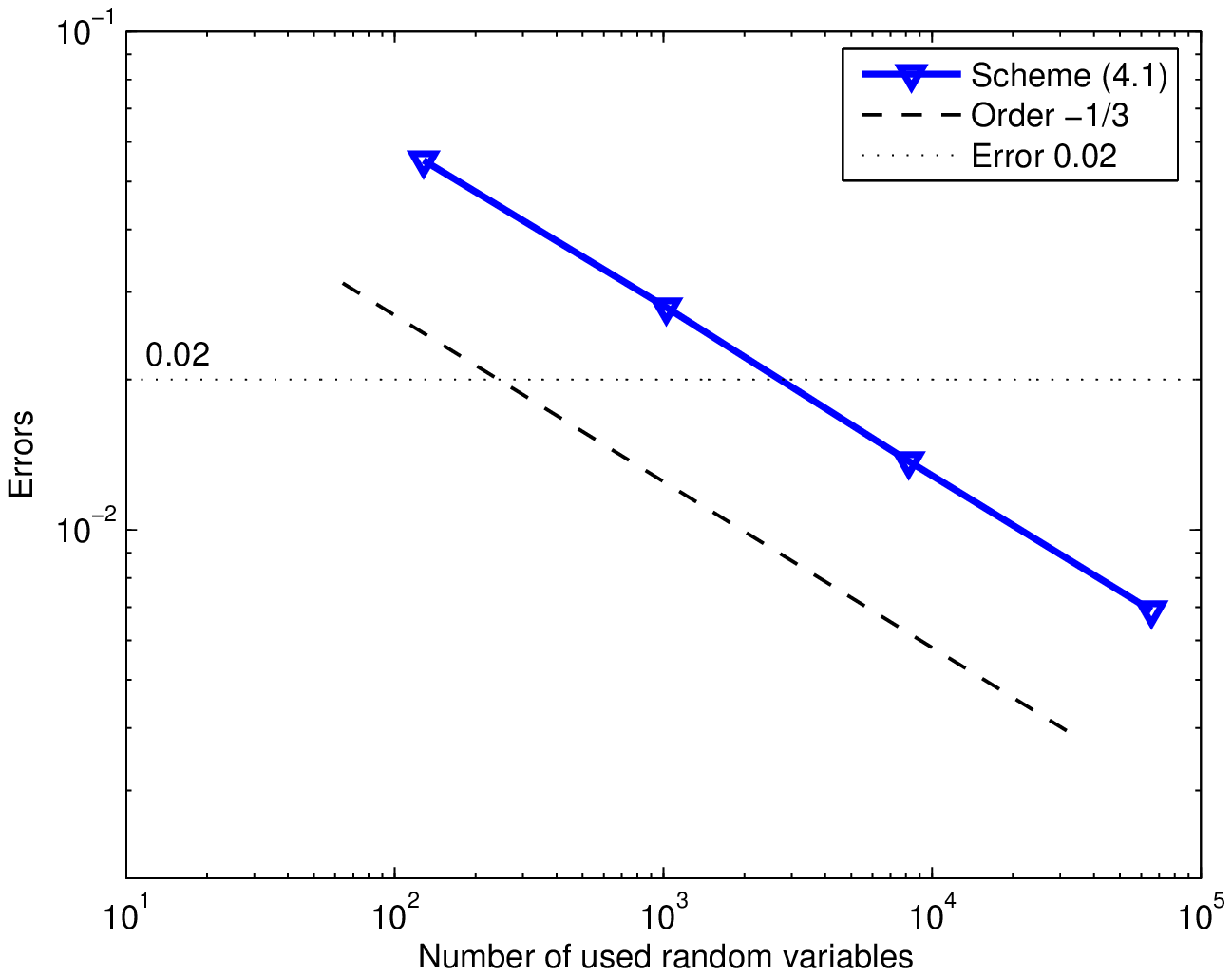}
         \includegraphics[width=2.5in,height=2.3in]
         {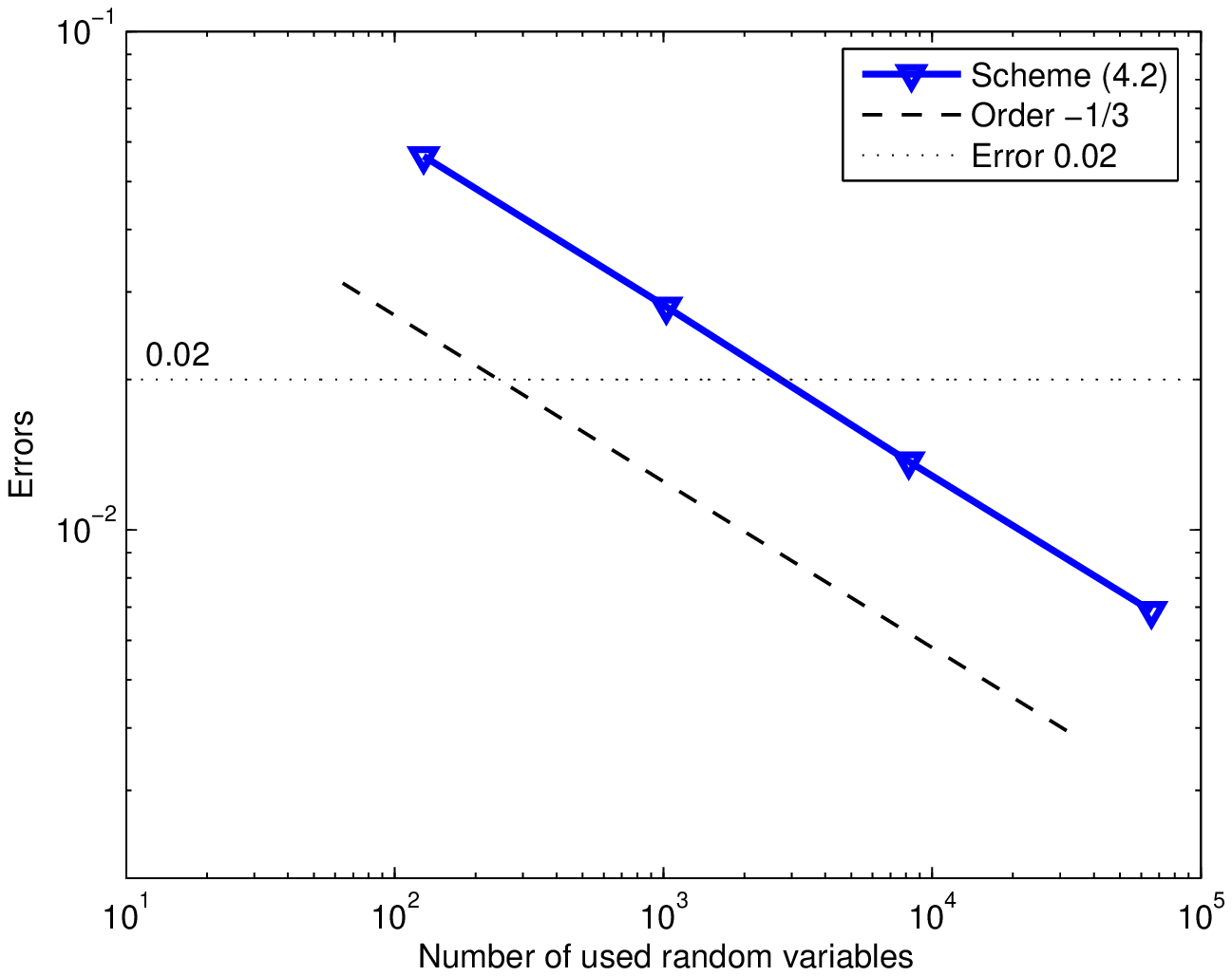}
         \caption{The overall approximation errors against number of used random variables.}
         \label{fig:overall.error}
\end{figure}

For the second example, we look at the nonlinear SWE \begin{equation}\label{eq:SWE.nonlinear2}
\left\{
    \begin{array}{lll}
    \frac{\partial^2 u}{\partial t^2} = \frac{\partial^2 u}{\partial x^2} +  \frac{1+u}{1+ u^2} + \dot{W}, \quad 0 < t \leq 1, \:\: x \in (0,1),
    \\
     u(0, x) = 0, \, \frac{\partial u}{\partial t}  (0,x) = 1, \: x \in (0,1),
     \\
     u(t, 0) = u(t,1) = 0,  \: t >0.
    \end{array}\right.
\end{equation}
Subsequently, we focus on the overall computational efforts of various fully discrete schemes, with the spectral Galerkin discretizations in space. We take the number of realizations of independent random variables needed for
approximations as a measure for the computational effort.
%To get the approximations $u_M^N$ via the numerical
%schemes mentioned above, one needs to generate
%$M\times N$ random variables.
Recall that the Crank-Nicolson-Maruyama (CNM) scheme and the stochastic trigonometric method (STM) converge with rate $\frac{1}{3} - \epsilon$ and order $\frac{1}{2} - \epsilon$ in time, respectively. In space, the spectral Galerkin method converges with order $\frac{1}{2} - \epsilon $ for arbitrarily small $\epsilon > 0$. In order to balance the errors in space and in time, we set $M = N^{\frac{3}{2}}$ for CNM scheme and $M = N$ for STM. Similarly, we set $M = N^{\frac{1}{2}}$ for the schemes \eqref{eq:full.discrete.concrete} and \eqref{eq:full.discrete.concrete2}. With these settings, the four schemes all result in an overall approximation error $O\big(N^{-\frac{1}{2} + \epsilon} \big)$.  The overall approximation errors produced by different schemes are listed in Table \ref{table:CNM}-\ref{table:EE}. Note that the ``exact'' solution is approximated by the spectral Galerkin method in space with $N_{exact} = 2^{12}$ and the STM scheme in time with $\tau_{exact} =  2^{-12}$. Again, 100 samples are used for the approximation of the expected values. In Fig.~\ref{fig:overall.error} we plot these approximation errors against the number of used random variables. The overall approximation errors of the four schemes all decrease at expected rates, as $N$ increases. For example, the overall computational errors of the schemes \eqref{eq:full.discrete.concrete} and \eqref{eq:full.discrete.concrete2} both decrease at  slope $-\frac{1}{3}$. This is an immediate consequence of our previous setting $M = N^{\frac{1}{2}}$. However, the overall errors of the other two methods exhibit decay rates of $-\frac{1}{5}$ and $-\frac{1}{4}$ as expected. Given a precision $\varepsilon = 0.02$, we are to compare the required computational costs for the above four schemes. One can first detect that the Crank-Nicolson-Maruyama scheme achieves the given precision  in the case $N = 2^8, M =2^{12}$ and thus requires to generate  $2^{20} = 1048576$ random variables. For the stochastic trigonometric method, $2^{18} = 262144$ ($N = M = 2^9$) random variables are needed to promise the precision. Our proposed schemes \eqref{eq:full.discrete.concrete} and \eqref{eq:full.discrete.concrete2}, however, achieve the given precision as $N = 2^8, M = 2^4$, which requires generation of only $2\times 2^{12} = 8192$ random variables. As discussed above, with the same precision, the proposed schemes \eqref{eq:full.discrete.concrete} and \eqref{eq:full.discrete.concrete2} can reduce the number of used random variables greatly and improve the computational efficiency significantly.

\begin{table}[htp]
\begin{center} \footnotesize
\caption{Computational errors of the Crank-Nicolson-Maruyama scheme with $M = N^{3/2}$}
\label{table:CNM}
\begin{tabular*}{8cm}{@{\extracolsep{\fill}}cccc}
\hline
$N = 2^2 $ & $ N = 2^4 $ & $N = 2^6$ & $N = 2^8$\\
\hline
0.13058 & 0.065411 & 0.032987 & 0.016622 \\
\hline
\end{tabular*}
\end{center}
\end{table}
\begin{table}[htp]
\begin{center} \footnotesize
\caption{Computational errors of the stochastic trigonometric method with $M =N$}
\label{table:STM}
\begin{tabular*}{8cm}{@{\extracolsep{\fill}}cccc}
\hline
$N = 2^6 $ & $ N = 2^7 $ & $N = 2^8$ & $N = 2^9$\\
\hline
0.054405 & 0.037954 & 0.026312 & 0.017867 \\
\hline
\end{tabular*}
\end{center}
\end{table}
\begin{table}[htp]
\begin{center} \footnotesize
\caption{Computational errors of the schemes \eqref{eq:full.discrete.concrete} and \eqref{eq:full.discrete.concrete2}, with $M = N^{1/2}$}
\label{table:EE}
\begin{tabular*}{10cm}{@{\extracolsep{\fill}}ccccc}
\hline
& $N= 2^4$& $ N = 2^6$ & $N = 2^8$ & N= $2^{10}$\\
 \hline
Scheme \eqref{eq:full.discrete.concrete} & 0.055098 & 0.027929 & 0.01372 & 0.0068861 \\
Scheme \eqref{eq:full.discrete.concrete2} & 0.05617 & 0.028007 & 0.013708 & 0.0068762 \\
\hline
\end{tabular*}
\end{center}
\end{table}

\section{Conclusion}
\label{sect:conclusions}

By incorporating suitable linear functionals of the noise process, two higher order fully discrete schemes have been devised in this work to
solve stochastic wave equations with additive space-time
white noise. Both theoretical convergence results and
numerical experiments show that the proposed schemes
can significantly reduce the computational costs,
compared with existing methods for the SWE \eqref{eq:SWE}.
%This is due to including more information on stochastic
%convolutions by exploiting linear functionals
%of noises in the presented schemes, which leads to a
%breakthrough of convergence rates in time.
%We mention that the new schemes can be also
%applied to SWE driven by a general
%additive noise, i.e., $Q$-Wiener process,
%on the condition that $\Lambda$ and $Q$ have a common
%eigenbasis. Recall that the stochastic trigonometric
%method incorporating only basic Wiener increments can
%achieve the rate $1 -$ in time in the case when the
%driven noise is trace class (i.e., $\mbox{Tr}(Q) <
%\infty$)\cite{CLS13}. Thus in that situation the
%strategy of using functionals of noises as discussed
%in this work is not necessary.
In the present work, we always assume that the drift coefficient satisfies a globally Lipschitz condition, which  was commonly used in the literature but excludes many important model equations in applications. In the future, we plan to address this issue and to investigate strong convergence
of numerical schemes for nonlinear SWEs in non-globally Lipschitz settings \cite{GI98,GI99}. In addition, another future direction of study is to construct higher order time integrators on the basis of the Taylor expansions of the solution process of SWEs \cite{Jentzen11}.

%\bibliographystyle{siam}
%\bibliography{../Bib/bibfile}

\end{document}